\documentclass[12pt]{article}

\pdfoutput=1

\usepackage{amsmath,amssymb,amsfonts,color}
\usepackage{graphicx,psfrag,epsf}
\usepackage{enumerate}
\usepackage[authoryear]{natbib}

\usepackage{url} 

\newcommand{\blind}{1}

\usepackage{geometry}
\def\sI{{\mathcal I}}                            
\def\mJ{{\mathcal J}}                            

\def\sP{{\mathcal P}}

\def\sC{{\mathcal C}}

\def\mJ{{\mathcal J}}

\def\vectorfontone{\bf}
\def\vectorfonttwo{\boldsymbol}
\def\vx{{\vectorfontone x}}                      
\def\vy{{\vectorfontone y}}                      
\def\vz{{\vectorfontone z}}                      %

\def\vone{{\vectorfontone 1}}
\def\vzero{{\vectorfontone 0}}

\def\vbeta{{\vectorfonttwo \beta}}               
\def\vgamma{{\vectorfonttwo \gamma}}             %
\def\vtheta{{\vectorfonttwo \theta}}             
\def\vmu{{\vectorfonttwo \mu}}                   
\def\matrixfontone{\bf}

\def\mI{{\matrixfontone I}}                      
\def\mJ{{\matrixfontone J}}                      %
\def\mP{{\matrixfontone P}}                      %
\def\mX{{\matrixfontone X}}                      
\def\bE{{\mathbb E}}                             
\def\bP{{\mathbb P}}                             
\def\bI{{\mathbb I}}                             

\def\ds{\displaystyle}

\def\argmax{\operatornamewithlimits{\text{argmax}}}

\definecolor{orange}{rgb}{1,0.5,0}
\definecolor{purple}{rgb}{0.75,0,1}
\definecolor{darkgreen}{rgb}{0,0.5,0}

\begin{document}

\def\spacingset#1{\renewcommand{\baselinestretch}%
	{#1}\small\normalsize} \spacingset{1}

\if1\blind
{
\centerline{\bf \Large Bayesian hypothesis tests with diffuse priors: }
\medskip
\centerline{\bf \Large  Can we have our  cake and eat it too?\footnote{ 
		The authors gratefully acknowledge funding from the \textit{Australian Research Council -- DECRA Fellowship  DE130101670}}}

\medskip
\medskip
	\centerline{\sc by John T. Ormerod${}^{(1,2)}$, Michael Stewart${}^{(1)}$}
	\medskip
	\centerline{\sc  Weichang Yu${}^{(1)}$
	and  Sarah E. Romanes${}^{(1)}$}	
	\medskip
	\medskip
	\centerline{\it ${}^{(1)}$ School of Mathematics and Statistics, University of Sydney, Sydney 2006, Australia}
	\medskip
	\centerline{\it ${}^{(2)}$ ARC Centre of Excellence for Mathematical \& Statistical Frontiers}

	\medskip
	\centerline{25th of October 2017}
}

\fi

\if0\blind
{
	\bigskip
	\bigskip
	\bigskip
	\begin{center}
		{\LARGE\bf Bayesian hypothesis tests with diffuse priors: 
			Can we have our  cake and eat it too?}
	\end{center}
	\medskip
} \fi

\bigskip
\begin{abstract}
We introduce a new class of priors for  Bayesian hypothesis testing, which we name ``cake priors''.
These priors circumvent Bartlett's paradox (also called the Jeffreys-Lindley paradox); the 
problem associated with the use of diffuse priors leading to nonsensical statistical inferences.
Cake priors allow the use of diffuse priors (having one's cake) while achieving theoretically 
justified inferences (eating it too). We demonstrate this methodology for Bayesian hypotheses tests for
scenarios under which the one and two sample $t$-tests, and linear models
are typically derived. 
The resulting Bayesian test statistic takes the form of a penalized 
likelihood ratio test statistic. By considering the sampling distribution under the null and 
alternative hypotheses we show for independent identically distributed regular parametric models that Bayesian hypothesis tests
using cake priors
are Chernoff-consistent, i.e., achieve zero type I and II errors asymptotically. 
Lindley's paradox is also discussed. We argue
that a true Lindley's paradox  will only occur with small 
probability for large sample sizes. 
\end{abstract}

\noindent%
Chernoff-consistent; cake priors; linear models;
asymptotic properties of hypothesis tests.
\vfill

\newpage
\spacingset{1.45} 

\section{Introduction}

Determining appropriate parameter prior distributions is of paramount importance in Bayesian 
hypothesis testing. Bayesian hypothesis testing often centres around the concept of a Bayes factor,
which was initially developed by \cite{Jeffreys1935,Jeffreys1961}, and later popularized by 
\cite{Kass1995}. The Bayes factor is simply the odds of the marginal likelihoods between two 
hypotheses and is analogous  to the likelihood ratio statistic in classical statistics,
where instead 
of maximizing the likelihoods with respect to the model parameters, the model parameters are 
marginalized out. In classical statistical theory testing a simple point null 
hypothesis against a composite alternative is routine. However, such hypothesis tests can pose severe 
difficulties in the Bayesian inferential paradigm where the Bayes factors may exhibit undesirable properties unless parameter prior 
distributions are chosen with exquisite care. This paper offers a solution to this difficulty.

Prior distributions can be chosen in an informative or uninformative fashion. Employing informative priors (either 
based on data from previously conducted experiments, or eliciting priors from subject matter 
experts) can be impractical, particularly when the number of parameters in the model is large. 
Furthermore, informative priors can be criticised on grounds that such priors are inherently 
subjective or may not let the data from the current experiment speak for itself. 
However, using alternative priors can also lead to problems.

One such problem occurs when using overly diffuse or flat improper priors.
In the former case, as priors become more diffuse the hypothesis corresponding to
the smaller model becomes increasingly favoured regardless of the evidence
provided by the data. This problem occurs due to the normalizing constants of the priors dominating
the expression for the Bayes factor, and is sometimes referred to
as Bartlett's paradox \citep[e.g,][]{Liang2008}
or the Jeffreys-Lindley paradox \citep[e.g,][]{Robert1993,Robert2014}, named after
the pioneering work of \cite{Jeffreys1935,Jeffreys1961}, \cite{Lindley1957}, and \cite{Bartlett1957}, whose
authors identified this and other related problems associated
with Bayes factors. Discussions of this paradox and the related
Lindley's paradox 
can be found in \cite{Aitkin1991},
\cite{Bernardo1999}, \cite{Sprenger2013}, \cite{Spanos2013}, and
\cite{Robert2014}.

An extension of Bartlett's paradox occurs in the limit where diffuse priors become
flat to the point of being improper. The use of flat improper priors gives rise to
arbitrary constants in the numerator and denominator of the Bayes factor
\citep[see ][]{DeGroot1973}. Such
arbitrary constants are problematic since, without suitable modification, they could 
be chosen by the analyst to suit any
preconceived conclusions preferred, and as such, are not suitable
for scientific purposes.
Techniques for selecting the arbitrary constants in Bayes factors in an acceptable
way when employing flat improper priors have been developed in several papers.
\cite{Bernardo1980}
proposes to derive a reference prior for the null hypothesis
by maximizing a measure of missing information.
\cite{Spiegelhalter1982},
and \cite{Pettit1992} use an imaginary
data device leading to the arbitrary constants cancelling with other
terms in the Bayes factor. A further approach to the problem of using diffuse priors was proposed by
\cite{Robert1993} who advocated for reweighing the prior odds to balance against parameter priors as prior hyperparameters become diffuse.

\cite{OHagan1995} considers the problem of using
flat improper priors in the calculation of the Bayes factor by splitting
the data into a training and testing set. The training set is used to
construct an informative prior to be used to calculate the Bayes factor
using the remaining portion of the data. These ideas have been refined
in \cite{OHagan1997},
\cite{Berger1996}, and \cite{Berger2001}.
A computational drawback of some of these approaches
is that the same model is fit multiple times. For models where Bayesian inferential procedures
are considered too slow for fitting a single model these approaches to Bayesian testing become
infeasible from a practical viewpoint.

Other Bayesian hypothesis testing approaches abandon the Bayes factor altogether
by constructing hypothesis testing criterion which only enter the criterion
through parameter posterior distributions themselves. These include information
criteria type approaches such as the Bayesian information
criterion (BIC) and the deviance information criterion (DIC).
The BIC or Schwarz's criterion uses a Laplace
approximation
where the prior term is assumed to be asymptotically negligible as the sample size grows
\citep{Schwarz1978}. 
The DIC involves
a linear combination of the log-likelihood evaluated at a suitably chosen Bayesian point
estimator and the posterior expectation of the log-likelihood
\citep{Spiegelhalter2002}. Under such a construction the DIC is not
dominated by the prior as prior hyperparameters diverge.
Similarly, posterior Bayes factors proposed by \cite{Aitkin1991}
are based on the posterior expectation of the likelihood function, rather than the joint likelihood
(comprising of the model likelihood and prior). Since  this 
only involves the prior in the calculation of the posterior  distribution, the prior does not dominate posterior Bayes factors. 
\cite{Berger1996} criticized this approach because it employs a
double use of the data that is not consistent with typical Bayesian logic.

An interesting alternative approach to Bayesian hypothesis testing is that suggested in Section 6.3 of
\cite{Gelman2013} who discuss examining
the posterior distribution of carefully chosen test statistics
such that
large values of a given test statistic provides
evidence against the null hypotheses. This idea is explored
more formally in \cite{GelmanEtal1996} and give rise to the 
concept of posterior predictive $p$-values, the probability
that a test statistic of posterior predictive values is greater
than the observed value of the test statistic.

Bayes factors in the context of linear model selection 
\citep{
Zellner1980b,
Mitchell1988,
George1993,
Fernandez2001,
Liang2008,
Maruyama2011,
Bayarri2012}
and generalized linear model selection
\citep{Chen2003,
	Hansen2001,
Wang2007,
Chen2008,
Gupta2009,
Bove2011,
Hanson2014,
Li2015}
have received an 
enormous amount of attention. While we defer discussion of the types of priors used in these 
contexts to Section \ref{sec:Example2} 
 we will draw special attention to \cite{Liang2008}.
\cite{Liang2008} considers several prior structures in the context of linear models.
They employ Zellner's $g$-prior \citep{Zellner1980b,Zellner1986} for the regression coefficients where $g$ is a prior
hyperparameter. They consider several choices for choosing $g$ including setting $g$ to various 
constants, selecting $g$ using a local and global empirical Bayes procedure, and
via placing a hyperprior on $g$. Their results suggest that in order for the
resulting Bayes factors to be well behaved (including model selection consistent) a hyperprior needs 
to be placed on $g$.


In this paper we will construct a new class of priors which was inspired by the
priors used in the context of linear and generalized linear models.
This class of priors is constructed in such a way as to mimic Jeffreys priors
\citep[which have the
desirable property that they are invariant under parameter transformations][]{Jeffreys1946} in the limit as a prior hyperparameter $g$ diverges.
In order to circumvent a Bartlett's like paradox from occurring, the rate
at which $g$ diverges is different in the null and alternative hypothesis
in such a way  that results in the cancellation of problematic terms 
in both the numerator and denominator of the
Bayes factor. 

Bayes factors using cake priors have several desirable properties. In the examples we 
consider the Bayes factor can be expressed as a difference in BIC values,
i.e., a penalized version of the likelihood ratio test (LRT) statistic.
Using properties of the LRT statistic we show that Bayesian
hypothesis tests are Chernoff-consistent in the sense of 
\cite[][Section 2.13]{Shao2003}, i.e., they achieve asymptotically
zero type I and type II errors as the sample size diverges. In contrast classical
hypothesis testing procedures are usually chosen to have a fixed type I error and are consequently not Chernoff-consistent. In this respect our Bayesian hypothesis 
tests are superior to classical procedures whose type I error is held fixed. 
Due to the above properties we call the priors we develop ``cake priors'' since they allow
the use of diffuse priors (having ones cake) while being able to perform
sensible statistical inferences (eating it too).
We will also discuss Lindley's paradox in the context of cake priors and argue that generally  
Lindley's paradox will only occur with vanishingly small probability
for large samples.

 In Section \ref{sec:Probelms}
we reintroduce Bayes factors, including the interpretation of Bayes
factors. In Section \ref{sec:LindlyBardletParadox} we discuss more
specifically the problems associated with Bayes factors,
including both Lindley's and Bartlett's paradoxes. In Section
\ref{sec:Cake} we describe cake priors and illustrate their use
in the context of
one sample tests for equal means (with unknown variance),
two sample tests for equal means (assuming unequal variances),
linear models, and one sample tests for equal
means (with known variance). In Section \ref{sec:theory} we derive
some asymptotic theory for our proposed of Bayesian hypothesis
tests. In Section \ref{sec:ArbitraryConstants} discuss
the relationship between cake priors and improper priors
and discuss how arbitrary constants can be introduced into the 
Bayes factor.
In Section \ref{sec:interpretation} we take a closer
look at the interpretation of Bayes factors in light of our findings.
In Section \ref{sec:discussion} we conclude.

\section{Bayes factors}
\label{sec:Probelms}

Bayes factors are a key concept in Bayesian hypothesis testing introduced by 
\cite{Jeffreys1935,Jeffreys1961}, although a similar concept was
also developed independently by \cite{Good1952}.
Suppose that we have
observed the data vector $\vx = (x_1,\ldots,x_n)^T$
which are observed samples from $\sP = \{ \, p_i( \, \cdot \, ) \, \colon \, i=1,\ldots,n \, \}$
and we have two hypotheses $H_0$ and $H_1$ representing two models
$\sP_j = \{ \, p_{ij}( \, \cdot \, | \vtheta_j, H_j) \, \colon \, i=1,\ldots,n \, \}$, $j=0,1$,
describing two potential distributions from which $\vx$ was drawn, i.e.,
\begin{equation}\label{eq:hypothesesGeneral}
H_0 \colon \sP  \in \sP_0
\qquad \mbox{versus} \qquad
H_1 \colon \sP  \in \sP_1.
\end{equation}

\noindent The models could potentially have distinct parameters from two distinct models and the models need not be nested.
Let $p(\vtheta_j|H_j)$ be the prior distribution under hypothesis $H_j$
for $j=0,1$. The Bayes factor is then defined as
$$
\mbox{BF}_{01} = \frac{p(\vx|H_0)}{p(\vx|H_1)} =
\frac{\int p(\vx|\vtheta_0,H_0)p(\vtheta_0|H_0) d\vtheta_0}{\int p(\vx|\vtheta_1,H_1)p(\vtheta_1|H_1) d\vtheta_1},
$$

\noindent where integrals are replaced with combinatorial sums for discrete
random variables. The posterior odds of $H_0$ to $H_1$ is defined by
$\mbox{PO}_{01} =  \mbox{BF}_{01} \times  p(H_0) /p(H_1),$
where the factor $p(H_0)/p(H_1)$ is the prior odds.  
Assuming $p(H_0) = p(H_1) = 1/2$, the Bayes factors
have the interpretation that
when $\mbox{BF}_{01}$  is   above 1 the hypothesis $H_0$ is favoured and
when $\mbox{BF}_{01}$ is below 1 the hypothesis $H_1$ is favoured.
However, if the prior odds is not equal to one then the posterior
odds should be the focus for inference.

A Bayesian hypothesis test function $T(\vx) \in \{ 0,1 \}$ is based on
$T(\vx) = I(\lambda_{\mbox{\scriptsize Bayes}} > 0)$ where
$\lambda_{\mbox{\scriptsize Bayes}} = - 2\ln\mbox{BF}_{01}$
(which we will call the Bayesian test statistic, analogous to the LRT statistic),
and
$$
T(\vx) = \left\{
\begin{array}{rl}
1 & \ds \mbox{implies $H_1$ is preferred; \quad and} \\
0 & \ds \mbox{implies $H_0$ is preferred}.
\end{array}
\right.
$$

As indicated in the above equation, the interpretation
of results based on Bayesian hypothesis tests is different
from the interpretation of frequentist tests.
Frequentist hypothesis testing which asks whether the data could
have plausibly been drawn from the null model (based on a chosen test statistic), without reference
	to an alternative model (although LRT statistics
	are, for example, {\em derived} with reference to a specific alternative model).
Furthermore, in the Bayesian paradigm preference towards a particular hypothesis is stated,
	rather than rejection of the null. Note that preference should not be confused with 
	endorsement.
	One can have a preference between two poorly fitting models without stating that
	either model fits the data well.


\cite{Kass1995}
offer an
interpretation of
$\lambda_{\mbox{\scriptsize Bayes}}$,
and $\mbox{BF}_{10} = 1/\mbox{BF}_{01}$
in Table \ref{tab:bayesFactorInterpretation}
in terms of strength of evidence against the null hypothesis. 
In Section \ref{sec:interpretation} we will take a closer
look at the interpretation of Bayes factors in light of the analysis in the current paper.
For the examples we consider, using the cake priors
described later,
the quantity $\lambda_{\mbox{\scriptsize Bayes}}$ will turn out
to be a penalized version of $\lambda_{\mbox{\scriptsize LRT}} = -2[\ell_0(\widehat{\vtheta}_0) - \ell_1(\widehat{\vtheta}_1)]$ 
(the LRT statistic for the hypotheses in (\ref{eq:hypothesesGeneral}) where $\ell_j(\vtheta_j)= \ln p(\vx|\vtheta_j,H_j)$ 
and the $\widehat{\vtheta}_j$'s are
the MLEs under $H_j$) given by
\begin{equation}\label{eq:BayesLRT}
\ds \lambda_{\mbox{\scriptsize Bayes}} = \lambda_{\mbox{\scriptsize LRT}} - \nu\ln(n)
\qquad \mbox{or} \qquad 
\lambda_{\mbox{\scriptsize Bayes}} = \lambda_{\mbox{\scriptsize LRT}} - \nu\ln(n) + O(n^{-1})
\end{equation} 

\noindent depending on the example, 
where $\nu$ is the difference in the number of parameters
in $H_0$ and $H_1$.
Intuitively one might expect $\lambda_{\mbox{\scriptsize Bayes}}$
and $\lambda_{\mbox{\scriptsize LRT}}$ to be related since the focus
of both approaches are based on the ratio of likelihoods, albeit
different likelihoods.

\begin{table}[ht]
	\begin{center}
		\begin{tabular}{c|c|c}
			$\lambda_{\mbox{\scriptsize Bayes}}$	& $\mbox{BF}_{10}$ & Strength of evidence  \\
			\hline
			$0$ to $2$	&  $1$ to $3$  & not worth more than a bare mention
			\\
			$2$ to $6$	&  $3$ to $20$ & positive
			\\
			$6$ to $10$	&  $20$ to $150$ & strong
			\\
			$>10$	& $>150$ & very strong
		\end{tabular}
	\end{center}
	\caption{Table of interpretation of Bayes factors offered by \cite{Kass1995}.}
	\label{tab:bayesFactorInterpretation}
\end{table}

\section{Paradoxes in Bayesian hypothesis testing}
\label{sec:LindlyBardletParadox}

Problems with Bayesian hypothesis testing based on Bayes factors, for particular combinations
of hypotheses and priors, have been identified as early as 1935 by \cite{Jeffreys1935}, and later by \cite{Lindley1957}, and
\cite{Bartlett1957}. 
As we will see for particular hypothesis
tests, when parameter priors are not chosen with care, the conclusions based on
Bayes factors will not be sensible.
To give some context for the ensuing discussion we will now consider the
 hypothesis testing problem introduced by \cite{Lindley1957}  in order
to illustrate potential problems.

\medskip
\noindent
{\bf Lindley's example:}
Consider the hypothesis test where the sample is modelled via
$x_i|\mu \sim N(\mu,\sigma^2)$,
$1\le i\le n$ independently,
where $\mu$ and $\sigma^2$ are the mean and variance
parameters respectively. Here $\mu$ is an unknown value to
be estimated and $\sigma^2$ is a fixed known constant.
Suppose that we wish to perform the hypothesis test
\begin{equation}\label{eq:hypZtest}
H_0\colon \mu = \mu_0
\qquad \mbox{versus} \qquad
H_1\colon \mu \ne \mu_0,
\end{equation}

\noindent where $\mu_0$ is a known constant.
Under $H_0$ the values of all model parameters are fixed
(so that under $H_0$ the model has zero unknown parameters), i.e., $H_0$
is a simple point null hypothesis.
Suppose
that for $H_1$  we employ the prior
$\mu|H_1 \sim N(\mu_0,\tau^2)$
%
where the prior variance $\tau^2$ is a known constant.  
The Bayes factor with the stated prior on $\mu$
is
\begin{equation}\label{eq:Lindley}
\mbox{BF}_{01} =
\sqrt{1 + \frac{n\tau^2}{\sigma^2}}
\exp\left[ 
- \frac{n z(\vx)^2}{2(n + \sigma^2/\tau^2)} 
\right],
\end{equation}

\noindent where $z(\vx) = \sqrt{n}(\overline{x} - \mu_0)/\sigma$ is
the standard $z$-test statistic \citep[see][]{Bernardo1999}. 
The $p$-value for this test is $\bP( \chi_1^2 > z(\vx)^2 )$.
\vspace{-6mm}\begin{flushright}$\square$\end{flushright}

\vspace{-2mm}
\noindent
If we were to choose $\mu|H_1$ as above then 
\cite{Lindley1957}  identified the following problem.

\begin{itemize}
	\item {\bf Problem I}: For any fixed $p$-value as $n\to\infty$ we have $\mbox{BF}_{01}\to \infty$.
	
\end{itemize}

\noindent 
Suppose that the observed value of $z(\vx)$ is large so that, for any reasonably chosen level $\alpha$,
the typical frequentist approach would reject the null hypothesis.
For this value of $z(\vx)$ a Bayesian procedure based on the above Bayes 
factor would prefer the null hypothesis for a sufficiently large $n$,
drawing a contradiction between the two inferential paradigms. 
Problem I is associated with Lindley's paradox,  
also referred to as the Lindley-Bartlett and the Jeffreys-Lindey paradox). For discussion of this see \cite{Smith1980,Aitkin1991,Bernardo1999,Sprenger2013,Robert2014}.

We now consider a second example posed by \cite{Sprenger2013}.
%
%
%
%

\medskip 
\noindent 
{\bf Sprenger's example:} \cite{Jahn1987}  used electronic and
quantum-mechanical random event generators with visual
feedback; the subject with alleged psychokinetic ability tries to
``influence'' the generator.
The number of ``successes'' was $s=52,263,470$ and
the number of trials was $n=104,490,000$.
Assuming independence of each trial we have
$x_i|\rho \sim \mbox{Bernoulli}(\rho)$, $1\le i\le n$ with $p\in[0,1]$.
If we test
\begin{equation}\label{eq:equalProportion}
H_0\colon \rho=0.5 \qquad \mbox{versus} \qquad H_1\colon \rho\ne 0.5,
\end{equation}

\noindent a rejection of $H_0$ leads to evidence that the
subject has psychokinetic ability. Using the data a classical hypothesis
testing approach leads to a $p$-value approximately equal to $0.0003$, leading to a rejection of the
null hypothesis for the $\alpha=0.05$ cut-off. A 95\% confidence interval for $\rho$ is $(0.50008, 0.50027)$.
A standard Bayesian hypothesis test uses the prior
$\rho\sim\mbox{Beta}(1/2,1/2)$ (the Jeffreys prior)
leads to:
\begin{equation}\label{eq:BayesBinomTest}
\lambda_{\mbox{\scriptsize Bayes}} =
2\ln\mbox{Beta}(1/2 + s,1/2 + n - s)
- 2\ln(\pi) + 2n\ln(2),
\end{equation}

\noindent where $s= \sum_{i=1}^n x_i$.  For the Bayesian test
$\lambda_{\mbox{\scriptsize Bayes}} \approx -5.86$, which implies
the null model is preferred and an apparent contradiction 
between inferential paradigms.

\subsection{Resolving Lindley's paradox}

We will now resolve Lindley's paradox in both of the above examples.

\medskip 
\noindent 
{\bf Resolving Lindley's example:}
We argue that Problem I for Lindley's example only occurs because it is assumed that the
$p$-value is held fixed,
and that a true Lindley's paradox only occurs with vanishingly small probability as $n\to\infty$. 
The $p$-value cannot be held fixed as $n\to\infty$ as its behaviour depends on the
data generation process. 
Let $\mX = (X_1,\ldots,X_n)^T$ be a random sample.
Consider the value of $\lambda_{\mbox{\scriptsize Bayes}}$ for Lindley's example 
as a function of this random sample, i.e., where
\begin{equation}
\label{eq:LambdaBayesExample1}
\ds \lambda_{\mbox{\scriptsize Bayes}}(\mX)
= \frac{n z(\mX)^2}{n + \sigma^2/\tau^2} - \ln\left( 1 + \frac{n\tau^2}{\sigma^2}\right).
\end{equation}

\noindent The first term 
on the right-hand side of (\ref{eq:LambdaBayesExample1}) depends on the data
generating process for $\mX$, whereas 
the second term is $O(\ln(n))$. Consider the two cases:
\begin{enumerate}
	\item If $X_i \stackrel{\mbox{\scriptsize iid}}{\sim} N(\mu_0,\sigma^2)$, i.e., the data is generated from $H_0$. Then
	$z(\mX)^2 \sim \chi_1^2 = O_p(1)$ and the $O(\ln(n))$ term dominates.
	Hence, as $n\to\infty$ we have $\lambda_{\mbox{\scriptsize Bayes}}(\mX) \to -\infty$
	implying $\bP(
	T(\mX)=0)\to 1$, i.e., the null hypothesis is preferred.
	
	\item If $X_i \stackrel{\mbox{\scriptsize iid}}{\sim} N(\mu_1,\sigma^2)$ for some $\mu_1$ with $\mu_1\ne\mu_0$, i.e., the data is generated from $H_1$. 
	Then 
	$z(\mX)^2 \sim \chi_1'^{2}( n(\mu_1 - \mu_0)^2/\sigma^2 )$
	where $\chi_\nu'^2(\lambda)$ is the non-central chi-squared distribution with degrees of freedom $\nu$
	and non-centrality parameter $\lambda$.
	Then $z(\mX)^2 = O_p(n)$ dominating $O(\ln(n))$ term. Then as $n\to\infty$ we have $\lambda_{\mbox{\scriptsize Bayes}}(\mX) \to \infty$  
	implying $\bP(
	T(\mX)=1)\to 1$, i.e., the alternative hypothesis is preferred.
\end{enumerate}

\noindent Note that 1. implies that a test based on the above Bayes factor has vanishing type I error as $n\to \infty$ and 2. 
implies that the Bayesian test is consistent in the sense of \cite[][Section 3.3]{lehmann2004}. Combining 1. and 2. implies that the test is Chernoff consistent (see Section \ref{sec:theory} for a formal definition).

\medskip 
\noindent 
{\bf Resolving Sprenger's example:}
Using properties of the beta function, the gamma function, and
Stirling's approximation  leads to approximating 
(\ref{eq:BayesBinomTest})  by
$$
\lambda_{\mbox{\scriptsize Bayes}}(\vx) = \lambda_{\mbox{\tiny LRT}}(\vx) - \ln(n) - \ln(\pi/2) + O(n^{-1}),
$$

\noindent where
$\lambda_{\mbox{\tiny LRT}}$ is the LRT statistic corresponding to
the hypotheses (\ref{eq:equalProportion}).
Again, Lindley's paradox occurs here if we consider $\lambda_{\mbox{\tiny LRT}}(\vx)$
(or equivalently the $p$-value)
to be fixed. If $\lambda_{\mbox{\tiny LRT}}$ is held fixed and $n$ diverges
then the null hypothesis will be preferred in the limit.
Let $\mX = (X_1,\ldots,X_n)^T$ be a random sample.
We will later show (see Section \ref{sec:theory}) that
$$
\lambda_{\mbox{\tiny LRT}}(\mX) = \left\{ \begin{array}{rl}
O_p(1) & \mbox{if $H_0$ is true; \, and} \\
O_p(n) & \mbox{if $H_1$ is true.} \\
\end{array} \right.
$$

\begin{figure}[ht]
	\centering
	\includegraphics[width=0.75\textwidth]{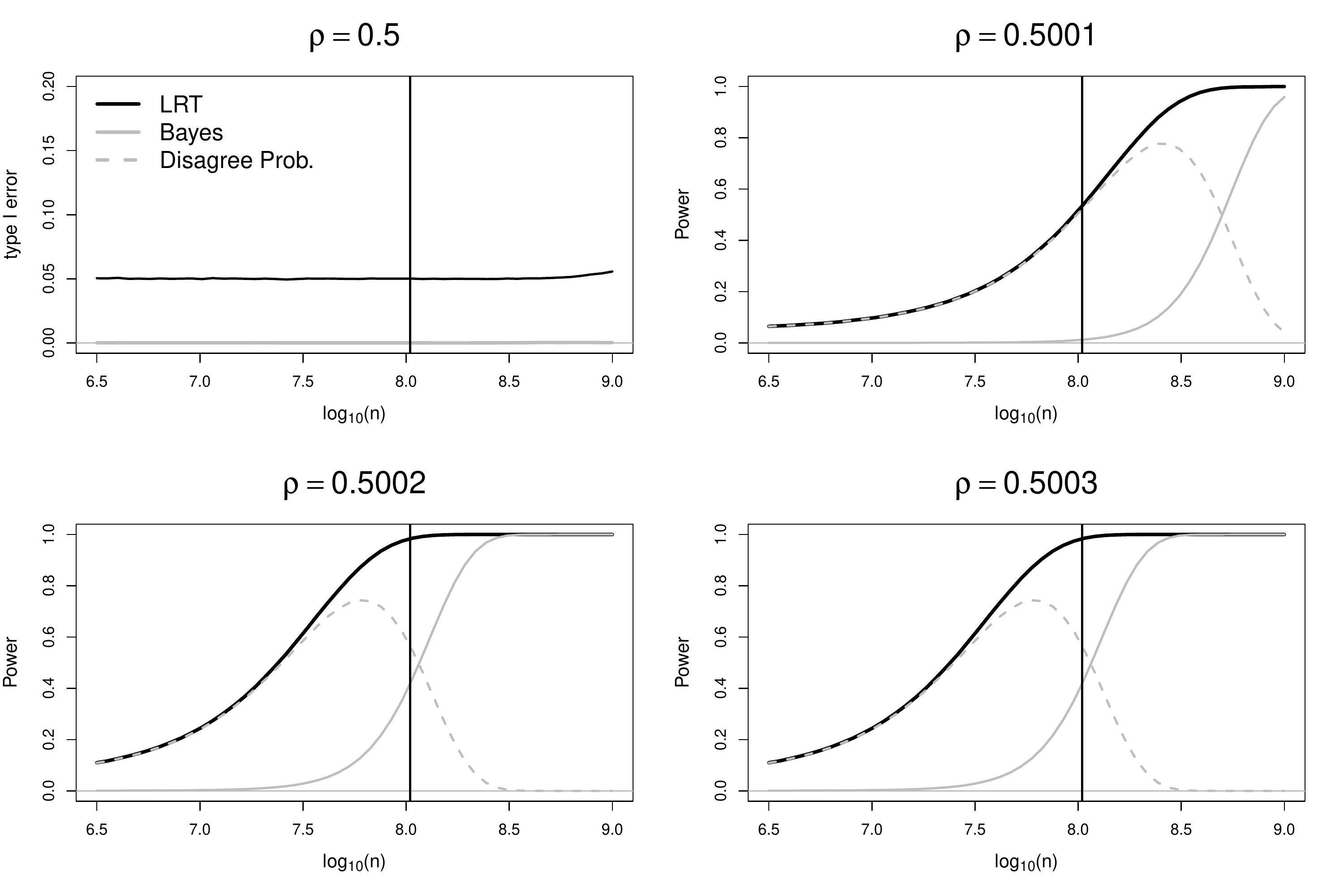}
	\caption{Empirical probabilities of rejecting the null hypothesis/preferring the alternative hypothesis for the simulation
		described in Section \ref{sec:LindlyBardletParadox} comparing the LRT
		and Bayesian tests when $\rho_{\mbox{\scriptsize true}} \in[0.5, 0.5004]$ and $n=104,490,000$.}
	\label{fig:01}
\end{figure}

\noindent Hence, as $n\to\infty$ we have $\lambda_{\mbox{\scriptsize Bayes}}(\mX) \to -\infty$
if $H_0$ is true, and $\lambda_{\mbox{\scriptsize Bayes}}(\mX) \to \infty$ if $H_1$ is true so that the test based on
$\lambda_{\mbox{\scriptsize Bayes}}$ is Chernoff consistent.
Figure \ref{fig:01} illustrates the empirical probabilities for rejecting the null
hypothesis (for the frequentist test at the 5\% level) or preferring the alternative hypothesis
(for the Bayesian test) based on simulating $10^6$ datasets with the true value of $\rho$ in the set
$\{ 0.5, 0.5001, 0.5002, 0.5003 \}$ for $n$ on a grid form $n=10^{6.5}$ to $n=10^9$.
The vertical line in Figure \ref{fig:01} indicates the actual value of $n$ in Sprenger's example 
and the dashed grey line illustrates the estimated empirical probability that the two tests disagree.
In Figure \ref{fig:01}
we see that when $H_0$ is false, as the sample size increases, both tests reject the null as
$n$ or $\rho$ grows. When $\rho=0.5$ the frequentist test accepts the null model at the 5\% level,
while the Bayesian test prefers the null model with very low probability.
Furthermore, at the actual value of $n$ in the experiment the disagreement between frequentist and Bayesian tests
could have occurred by chance with relatively high probability. However, for much larger
large $n$ such a disagreement will only
occur with low probability when $H_1$ is true.
\vspace{-8mm}\begin{flushright}$\square$\end{flushright}

We now note for Lindley's example that the  test based on $\lambda_{\mbox{\scriptsize Bayes}}$ is an LRT test with an extremely small level $\alpha$ given by
$\alpha = \bP[ \chi_1^2 > \{ 1 + \sigma^2/(n\tau^2)\} \ln(1 + n\tau^2/\sigma^2)]$
(when $\tau^2$ is large)
while for Sprenger's example  that the  test based on $\lambda_{\mbox{\scriptsize Bayes}}$ is in also an LRT test with asymptotic level 
$\alpha = \bP[ \chi_1^2 > \ln(n) + \ln(\pi/2) ]$. Hence, we note, 
as was also argued by \cite{Naaman2016} and noted by \cite{Lindley1957}, that Lindley's paradox can also
be resolved by letting the level of the test $\alpha$ in the frequentist test
tend to zero as $n\to\infty$. 
We discuss Lindley's paradox more generally in Section 
\ref{sec:comparingTwoModels}.

The above expression for the level of the test for Lindley's example draws attention to a second problem for Lindley's example
which does not occur in Sprenger's example.
\begin{itemize}
	\item {\bf Problem II}: For any given $p$-value
	as $\tau^2\to\infty$ we have $\mbox{BF}_{01}\to \infty$.
\end{itemize}

\noindent This problem occurs because as the prior variance increases the level of the test
decreases.

Problem II is referred to as
Bartlett's paradox in \cite{Liang2008} and the
Jeffreys-Lindley paradox in Robert (2014).
Bartlett's paradox is  paradoxical since
as $\tau^2$ becomes  large the prior on $\mu$ becomes increasingly vague regarding the location of $\mu$.
However, in the attempt to be vague about the location of $\mu$ the prior becomes ``informative'' in favouring $H_0$ as the preferred hypothesis, again,
regardless of the evidence provided by the data.
Unlike Lindley's paradox, we believe that Bartlett's paradox is a real problem in practice since the use of diffuse 
priors can sometimes lead to testing procedures with extremely small power. Our proposed cake priors described in the next section circumvent this problem.

\section{Cake priors}
\label{sec:Cake}

Consider the general hypotheses (\ref{eq:hypothesesGeneral}).
Let $d_0$  and $d_1$ be the dimensions of $\vtheta_0$ and
$\vtheta_1$  respectively.
For the time being we will assume that $0 < d_0 \le d_1$ (later we will consider $0 \le d_0 \le d_1$).
Define the observed information and Fisher information
matrices as
$\mJ(\vtheta) = -
\nabla_\vtheta^2 \ell(\vtheta)$ and 
$\sI(\vtheta) = \bE_{x_i|\vtheta}\left[ -
\nabla_\vtheta^2 \ln p(x_i|\vtheta) 
\right]$
respectively.
Define the mean observed information matrix as
$\widetilde{\mJ}(\vtheta) = n^{-1}  \mJ(\vtheta)$ and the
mean expected observed information matrix as
$\widetilde{\sI}(\vtheta) = n^{-1} \bE_{\vx|\vtheta}\left[ \mJ(\vtheta)
\right]$. 
We will denote the Fisher information matrix under the null
and alternative hypotheses as $\sI_0(\vtheta_0)$ and
$\sI_1(\vtheta_1)$ with similar use of subscripts to denote
similar quantities such as $\mJ$, $\widetilde{\sI}$ and
$\widetilde{\mJ}$.
We define a Jeffreys prior as any density for $\vtheta$ such that
$p(\vtheta) \propto |\sI(\vtheta)|^{1/2}$.

We construct
cake priors   using the following ingredients:
\begin{enumerate}
	\item Define the priors
	\begin{equation}\label{eq:cakePrior}
	p(\vtheta_j| H_j; g_j)= \exp\left[
	- \tfrac{d_j}{2}\ln(2\pi g_j)
	+ \tfrac{1}{2}\ln|\mP_j(\vtheta_j)|
	- \tfrac{1}{2g_j}\vtheta_j^T\mP_j(\vtheta_j)\vtheta_j
	\right],
	\end{equation}
	
	\noindent where $\mP_j(\vtheta_j)$ is a prior precision matrix (assumed to be full rank).
	For all of the examples considered in this paper
	we will use $\mP_j(\vtheta_j) = \widetilde{\sI}_j(\vtheta_j)$.

	\item Set $g_j = h^{1/d_j}$ where $h$ is a common hyperparameter.
	
	
	\item Calculate the Bayes Factor as
	$$
	\ds \mbox{BF}_{01}(h) = \left[ \ds
		\int p(\vx|\vtheta_0,H_0)p(\vtheta_0|H_0 ; h^{1/d_0})    d\vtheta_0
	\right] \Big/ \left[ \ds
		\int p(\vx|\vtheta_1,H_1)p(\vtheta_1|H_1 ; h^{1/d_1})  d\vtheta_1	
	\right].
	$$

	\item Optional: Let $h\to\infty$ if flat priors are desired.
\end{enumerate}

\noindent
When $\mP_j(\vtheta_j) \propto \sI_j(\vtheta_j)$, $j=0,1$,
(\ref{eq:cakePrior})
leads to a Bayes Factor, in the limit as $g_j\to\infty$, that would have been obtained if a Jeffrey's prior is used. 
When $\mP_j(\vtheta_j) \propto \sI_j(\vtheta_j)$, $j=0,1$,
(\ref{eq:cakePrior}) are Jeffreys priors in the limit as $g_j\to\infty$. 
Letting $g_j\to\infty$ would be problematic if not for 2. which leads
to certain terms involving $h$ cancelling in the Bayes factor.
As $h\to \infty$, the priors on $\vtheta_j$ are made diffuse, but at a rate 
that depends on the $d_j$'s.
We keep 4. optional due to the fact that particular Bayesian 
procedures may require
proper priors. 
For particular examples in the coming subsections the above steps will raise complications.
These include:
(A) Model parameters may not be defined on the whole real line, e.g., variances. Priors of the form (\ref{eq:cakePrior}) are not appropriate for such
	model parameters; 
(B) The priors $p(\vtheta_j|g_j,H_j)$ may not be proper densities; and
(C) If $d_0 = 0$  using $g = h^{1/d_j}$ is problematic.
 The examples in the subsections below will illustrate how each of these
complications can be handled. 
(A) \& (B) will be dealt with in sections \ref{sec:Example1}
and \ref{sec:Example5}. Complication (C) will be dealt with in
Section \ref{sec:zeroCase}.

We will now give some intuition for how cake priors avoid Bartlett's paradox via the
following heuristic argument. Let $\mP_j(\vtheta_j) \equiv \mP_j$, i.e., the prior precision matrices 
are constant, then letting
$\mbox{BF}_{01} = 
\left[ \int p(\vx|\vtheta_0,H_0)p(\vtheta_0|H_0 ; g_0)    d\vtheta_0
\right] / \left[ 
	\int p(\vx|\vtheta_1,H_1)p(\vtheta_1|H_1 ; g_1)  d\vtheta_1\right],$
(which depends on $g_0$ and $g_1$ rather than $h$). Then the Bayesian test statistic is
$$
\begin{array}{l}
\ds \lambda_{\mbox{\scriptsize Bayes}}\\
\quad = 
2\ln\left[ 
\frac{\ds
	\int \exp\left\{ \ell_1(\vtheta_1)  - \tfrac{\vtheta_1^T\mP_1\vtheta_1}{2g_1} \right\} d\vtheta_1		
}{\ds
	\int \exp\left\{ \ell_0(\vtheta_0)   - \tfrac{\vtheta_0^T\mP_0\vtheta_0}{2g_0} \right\} d\vtheta_0
}
\right] 
+ d_0\ln(2\pi g_0)
- d_1\ln(2\pi g_1) 
+ \ln(|\mP_1|/|\mP_0|)
\\ [4ex]
\quad \ds  
= 
2\ln\left[ 
\frac{
	\int p(\vx|\vtheta_1,H_1) d\vtheta_1
}{
	\int p(\vx|\vtheta_0,H_0) d\vtheta_0		
}
\right] 
+ d_0\ln(2\pi g_0)
- d_1\ln(2\pi g_1) 
+ \ln(|\mP_1|/|\mP_0|) + O(g_0^{-1} + g_1^{-1})
\end{array} 
$$

\noindent where the second line is obtained using a Taylor series argument in $g_0$ and $g_1$.
Ignoring the dependency of $O(g_j^{-1})$ terms on the $\vtheta_j$'s, using Laplace's method 
on the numerator and denominator of the first term in the second line above, 
and setting $\mP_j = \widetilde{\mJ}_j(\widehat{\vtheta}_j)$ (where $\widehat{\vtheta}_j$ are 
the MLEs for the $\vtheta_j$'s), leads to  
\begin{equation}\label{eq:Jeffreys}
\begin{array}{rl}
\ds \lambda_{\mbox{\scriptsize Bayes}}
& \ds = \lambda_{\mbox{\scriptsize LRT}} - \nu\ln(n)
+ d_0\ln(g_0) - d_1\ln(g_1) 
+ O(g_0^{-1} + g_1^{-1} + n^{-1}).
\end{array} 
\end{equation}

\noindent The $O(n^{-1})$ error follows from the relative error of the Laplace's method 
applied to the numerator and denominator \citep{Tierney1989,KassEtal1990}. 
Suppose that $g_0=g_1 =g$. Then using
the asymptotic $\chi_\nu^2$ distribution $\lambda_{\mbox{\scriptsize LRT}}$ the
level of the test using (\ref{eq:Jeffreys}) is
$\alpha = \bP[ \chi_\nu^2 \ge \nu\ln(n g)].$
So that again we see that the power of the test goes to 0 as $g\to\infty$.
Here also we see that setting $g_0$ to be a large constant (making the prior for $\vtheta_0$ 
diffuse) leads the test to preferring $H_0$ while making $g_1$ large leads to preferring
$H_1$. Hence, the relative rates that $g_0$ and $g_1$ diverge must be considered.

Setting $g_j= h^{1/d_j}$ means that $d_0\ln(g_0) = d_1\ln(g_1)$ and leads to
$\lambda_{\mbox{\scriptsize Bayes}} 
= \lambda_{\mbox{\scriptsize LRT}} - \nu\ln(n) + O(h^{-1/d_0} + h^{-1/d_1}  + n^{-1})$.
For sufficiently large $h$ and $n$ we have 
$\lambda_{\mbox{\scriptsize Bayes}} \approx \lambda_{\mbox{\scriptsize LRT}} - \nu\ln(n)$.
The level of the test becomes approximately $\alpha = \bP[ \chi_\nu^2 \ge \nu\ln(n)]$.
We recognise the above arguments are informal in nature and will shortly illustrate cake priors
in concrete examples. 

Lastly, well us briefly discuss the choice of $\mP_j$. Setting $\mP_j = \mI$ leads to 
$\lambda_{\mbox{\scriptsize Bayes}} 
= \lambda_{\mbox{\scriptsize LRT}} - \nu\ln(n) 
+ \ln|\widetilde{\mJ}_0(\widehat{\vtheta}_0)| 
- \ln|\widetilde{\mJ}_1(\widehat{\vtheta}_1)|
+ O(h^{-1/d_0} + h^{-1/d_1}  + n^{-1}).$
This would be undesirable because of the additional computational burden
of the log-determinant terms (which can be considerable in some contexts), and because if $\lambda_{\mbox{\scriptsize LRT}} \approx \nu\ln(n)$ we would
prefer the model with larger $\ln|\widetilde{\mJ}_j(\widehat{\vtheta}_j)|$, i.e.,
larger standard errors. For this reason we would like 
$\mP_j \approx \widetilde{\mJ}_j(\widehat{\vtheta}_j)$ so that at least
approximate cancellation occurs.


\subsection{One sample test for equal means (with unknown variance)}
\label{sec:Example1}

Consider the hypothesis test (\ref{eq:hypZtest})   where
$x_i|\mu,\sigma^2 \stackrel{\mbox{\scriptsize iid}}{\sim} N(\mu,\sigma^2)$,
$1\le i\le n$,
where $\mu$ and $\sigma^2$ are the mean and variance
parameters respectively. Suppose now that both
$\mu$ and $\sigma^2$ are unknown parameters to be estimated (unlike the example in Section \ref{sec:Probelms}
where $\sigma^2$ was assumed to be known).

The mean expected information matrices for the null
and alternative hypothesis respectively are
$$
\widetilde{\sI}_0(\sigma^2) = 1/(2\sigma^{4})
\qquad \mbox{and} \qquad
\widetilde{\sI}_1(\mu,\sigma^2) = \left[ \begin{array}{cc}
\sigma^{-2} & 0 \\
0 & 1/(2\sigma^{4})
\end{array}
\right],
$$

\noindent which coincide with the Fisher information matrices in this case.
Hence, Jeffreys priors for $\mu$ and $\sigma^2$ are
$p(\mu)\propto 1$ and  $p(\sigma^2)\propto (\sigma^2)^{-1} I(\sigma^2>0)$ respectively.

We cannot directly use the methodology outlined in Section \ref{sec:Cake}
as $\sigma^2>0$. To handle this complication
we use the transformation $\sigma^2 = \exp(s)$. Under this transformation 
the mean expected information matrices become:
$\widetilde{\sI}_0(s) = 1/2$ and 
$\widetilde{\sI}_1(\mu,s) = \mbox{diag}(\exp(-s),1/2)$.
Using the steps for Section \ref{sec:Cake} under this transformation we have $s|H_0 \sim N(0,2g_0)$, i.e., $p(s|H_0) = [2\pi(2g_0)]^{-1/2} \exp[-s^2/(4g_0)]$.
Transforming back to the $\sigma^2$ parametrisation
gives 
$ \sigma^2|H_0 \sim LN(0,2g_0),$ 
with density
which is a Jeffreys prior for $\sigma^2$ in the limit as $g_0\to \infty$.

Note that for $H_1$ the upper left entry of $\widetilde{\sI}_1(\mu,s)$ depends on $s$.
This implies a conditional dependence of $\mu$ on
$s$ which leads to $\mu|s,H_1 \sim N(0,g_1\exp(s))$ and
$s|H_1 \sim N(0,2g_1)$. Transforming back
from the parametrisation in $s$  to the parametrisation
using $\sigma^2$
gives
$\mu|\sigma^2,H_1 \sim N(0,g_1\sigma^2)
$ and $
\sigma^2|H_1 \sim LN(0,2g_1).$
 Note again that these are both Jeffreys priors in the
limit as $g_1\to \infty$.

The marginal distributions of
$\vx$ given $H_0$ and $H_1$ are 
$$
\begin{array}{ll}
\ds p(\vx|H_0)
= \int_0^\infty
\frac{1}{\sqrt{2\pi}^{n}}\exp\left[
- \tfrac{n+2}{2}\ln(\sigma^2)
- \tfrac{n\widehat{\sigma}_{0}^2}{2\sigma^2}
- \tfrac{\ln(4\pi g_0)}{2}
-\tfrac{(\ln \sigma^2)^2}{4g_0}
\right] d\sigma^2, \qquad \mbox{and} \\ [2ex]
\ds p(\vx|H_1)
= \int_0^\infty
\frac{1}{\sqrt{2\pi}^{n}}\exp\left[
- \tfrac{n+2}{2} \ln(\sigma^2)
- \tfrac{n\widehat{\sigma}_{g}^2}{2\sigma^2}
- \tfrac{\ln(4\pi g_1^2)}{2}
- \tfrac{(\ln \sigma^2)^2}{4g_1}
- \tfrac{1}{2}\ln\left(n + \tfrac{1}{g_1}\right)
\right] d\sigma^2,
\end{array}
$$

\noindent where $\widehat{\sigma}_0^2 = n^{-1}\|\vx - \mu_0\vone\|^2$, and
$\widehat{\sigma}_g^2 = n^{-1}[ \|\vx\|^2 - (n\overline{x})^2/(n + g^{-1})]$. Suppose that we were to use $g = g_0 = g_1$ and let $g\to \infty$.
Then we would see a manifestation of Bartlett's paradox where the null hypothesis
is favoured since $\mbox{BF}_{01} \to\infty$ as $g\to\infty$.
If we instead use $g_0 = h$ and $g_1 = h^{1/2}$
then the Bayes factor simplifies to
$$
\begin{array}{rl}
\mbox{BF}_{01}(h) =
& \ds
\frac{\ds \int_0^\infty
	\exp\left[ - \left( \tfrac{n}{2} + 1 \right)\ln(\sigma^2) - \frac{n\widehat{\sigma}_{0}^2}{2\sigma^2}
	-\frac{(\ln \sigma^2)^2}{4h}
	\right] d\sigma^2}{\ds \int_0^\infty
	\exp\left[ - \left( \tfrac{n}{2} + 1 \right)\ln(\sigma^2)
	- \frac{n\widehat{\sigma}_{h^{1/2}}^2}{2\sigma^2}
	- \frac{(\ln \sigma^2)^2}{4h^{1/2}}
	- \tfrac{1}{2}\ln(n + h^{-1/2})
	\right] d\sigma^2} 
\end{array}
$$

\noindent which can be evaluated using univariate quadrature
or other methods for any fixed $h>0$.

Since both the above integrands
are monotonic as a function of $h$ with a well defined limit
as $h\to\infty$ we  can apply the monotonic convergence theorem
and take the limit $h\to\infty$ inside both integrals.  After simplifications including $\widehat{\sigma}_h^2 \to \widehat{\sigma}_{1}^2$ where
$
\widehat{\sigma}_{1}^2 =
n^{-1}\|\vx - \overline{x}\vone\|^2$
is the MLE for $\sigma^2$ under the alternative hypothesis
we obtain
$\lambda_{\mbox{\scriptsize Bayes}}  =
\lambda_{\mbox{\scriptsize LRT}}
- \ln(n),
$
 where
$\lambda_{\mbox{\scriptsize LRT}} =
n\ln( \widehat{\sigma}_{0}^2 )
- n\ln( \widehat{\sigma}_{1}^2 )$
is the LRT statistic corresponding to the hypothesis
(\ref{eq:hypZtest}). 
%

As $h\to\infty$ the parameter posterior distributions are given by
$$
\sigma^2|\vx,H_0 \sim \mbox{IG}\left(\tfrac{n}{2},\tfrac{n}{2}\widehat{\sigma}_{0}^2\right),
\quad
\mu|\vx,H_1 \sim t_n(\overline{x},n^{-1}\widehat{\sigma}_{1}^2),
\quad \mbox{and}
\quad
\sigma^2|\vx,H_1 \sim \mbox{IG}\left(\tfrac{n}{2},\tfrac{n}{2}\widehat{\sigma}_{1}^2\right).
$$


As a computational short-cut if the integrand is
a monotonic function of $h$ with a well defined limit
as $h\to\infty$ we will write
$$
\begin{array}{l}
\ds p(\vx|H_0)
\ds \stackrel{h\Rightarrow\infty}{=}
 \int_0^\infty
\exp\left[
- \left(\tfrac{n}{2} + 1\right)\ln(\sigma^2)
- \tfrac{n\widehat{\sigma}_{0}^2}{2\sigma^2}
\right] d\sigma^2, \qquad \mbox{and} \\ [1ex]
\ds p(\vx|H_1)
\ds \stackrel{h\Rightarrow\infty}{=} \int_0^\infty
\exp\left[
- \left(\tfrac{n}{2} + 1\right)\ln(\sigma^2)
- \tfrac{n\widehat{\sigma}_{1}^2}{2\sigma^2}
- \tfrac{1}{2}\ln(n)
\right] d\sigma^2,
\end{array}
$$

\noindent where the notation $\stackrel{h\Rightarrow\infty}{=}$ is used to denote
``equality in the limit as $h \rightarrow \infty$ after terms related to $h$ cancel in the numerator and denominator in the Bayes factor,  or terms related to $h$
vanish as $h$ diverges in the Bayes factor.''
The above expressions can be more easily simplified
using standard results to reach the same expression for $\lambda_{\mbox{\scriptsize Bayes}}$ as above.

\begin{figure}[ht]
	\centering
	\includegraphics[width=0.95\textwidth]{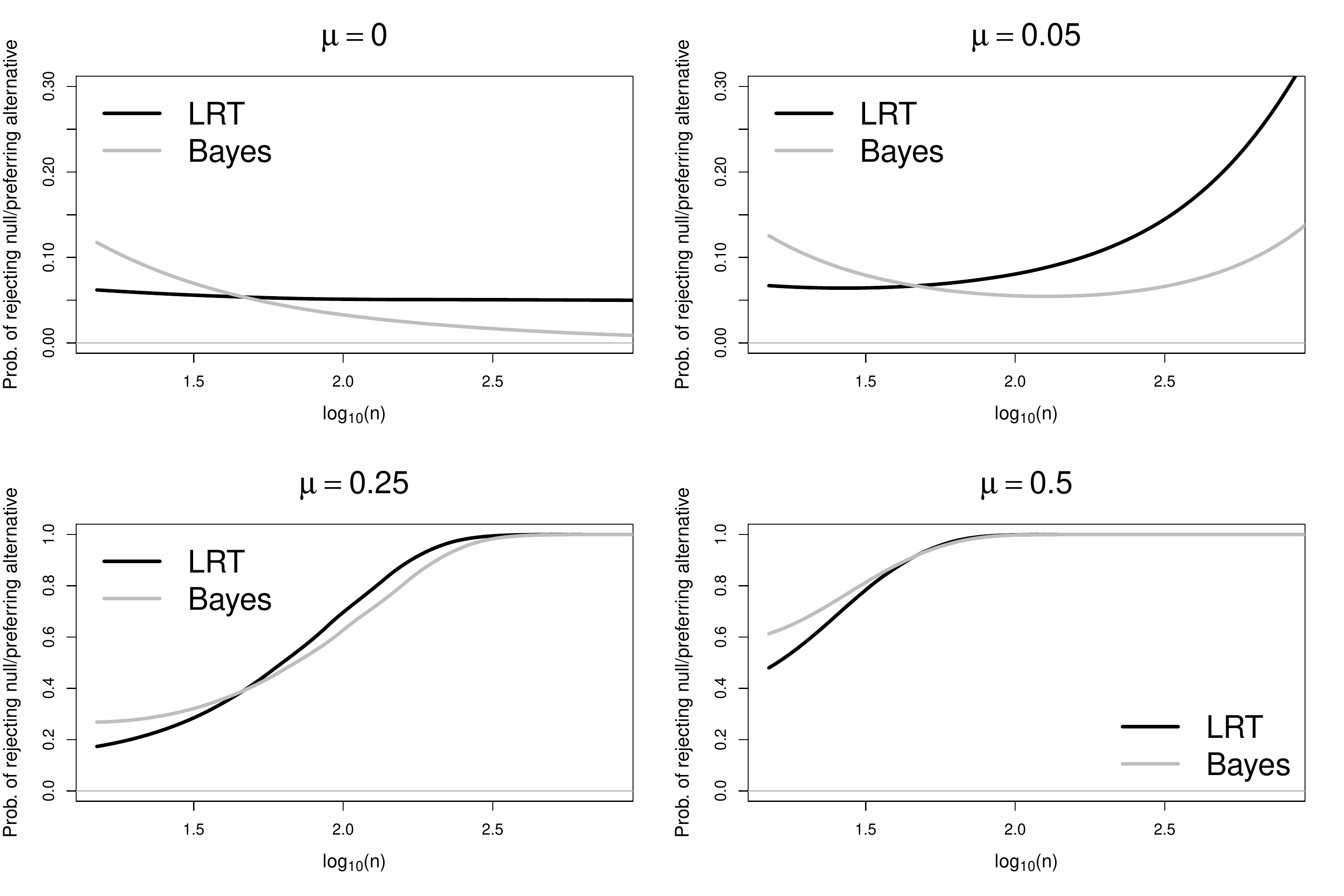}
	\caption{Empirical probabilities of rejecting the null hypothesis/preferring the alternative hypothesis for the simulation
	described in Section \ref{sec:Example1} comparing the LRT
	and Bayesian tests when $\mu_{\mbox{\scriptsize true}} \in\{ 0, 0.05, 0.25, 0.5\}$.}
	\label{fig:02}
\end{figure}

We conduct the following short simulation study to illustrate
the differences between the LRT and the Bayesian
test for this problem. Letting $\mu_0 = 0$ we simulate a single
set of data from
$x_i\sim N(\mu_{\mbox{\scriptsize true}},1)$, $1\le i\le n$.
After simulating $10^6$ such datasets for all values
of $\mu_{\mbox{\scriptsize true}}$ in the set $\{ 0, 0.05, 0.25, 0.5\}$
and a grid of $n$ from $n=15$ to $n=1000$ we plot in Figure \ref{fig:02}
the empirical
probabilities of rejecting the null hypothesis (for the LRT test)
using $\alpha = 0.05$ or preferring the alternative hypothesis
(for the Bayesian test).

From Figure \ref{fig:02} we see empirically that the type I error
of the Bayesian test is tending to 0 as $n$ grows when $H_0$ is true,
whereas the LRT test has, by design, a type I error of $0.05$.
When $H_1$ is true and $\ln(n) < \chi^2_{1,\alpha}$
the Bayesian test is more powerful than the LRT test, and when
$H_1$ is true and  $\ln(n) > \chi^2_{1,\alpha}$
the LRT test is more powerful than the Bayesian test.
When $\mu_{\mbox{\scriptsize true}} \in \{ 0.25, 0.5\}$
both tests appear to have power tending to 1 as $n$ grows.
Lastly, for the case $\mu_{\mbox{\scriptsize true}} = 0.5$
when $\ln(n) > \chi^2_{1,\alpha}$ both have very similar
power.

\subsection{Two sample test for equal means}
\label{sec:Example5}

\noindent
Suppose we have data $\vx = (x_1,\ldots,x_n)^T$. We want to test whether the first $n_0$
samples $\vx_0 = (x_1,\ldots,x_{n_0})^T$ from class 0
come from the same normal population as the second $n_1$ samples $\vx_1 = (x_{n_0+1},\ldots,x_{n})^T$ from class 1 with $n_0 + n_1 = n$.
We wish to test
\begin{equation}\label{eq:hypothesis2}
\begin{array}{l} 
H_0 \colon x_i|\mu,\sigma^2 \sim N(\mu,\sigma^2),
\ 1\le i\le n,
\qquad \mbox{versus} 
\\  
H_1 \colon
\left\{ \begin{array}{l}
x_{i}|\mu_0,\sigma_0^2 \sim N(\mu_0,\sigma_0^2), \ 1\le i\le n_0,  \\ [1ex]
x_{i}|\mu_1,\sigma_1^2 \sim N(\mu_1,\sigma_1^2), \  n_0+1\le i\le n,
\end{array}
\right.
\end{array}
\end{equation}

\noindent where $\mu$, $\sigma^2$, $\mu_0$, $\sigma_0^2$, $\mu_1$ and $\sigma_1^2$
are the means and variances under the one and two group hypotheses respectively.
Here $\vtheta_0 = (\mu,\sigma^2)^T$ with $d_0=2$
and $\vtheta_1 = (\mu_0,\mu_1,\sigma_0^2,\sigma_1^2)^T$ with $d_1=4$.
%
\noindent Using similar arguments as in Section \ref{sec:Example1}
with $\mP_0(\vtheta_0) = \widetilde{\sI}_0(\vtheta_0) 
= \mbox{diag}[ \sigma^{-2}, 1/(2\sigma^{4})]$ and
$\mP_1(\vtheta_1) = 
\widetilde{\sI}_1(\vtheta_1) = \mbox{diag}[
n_0/(n\sigma_0^{2}),
n_1/(n\sigma_1^{2}),
n_0/(2n\sigma_0^{4}),
n_1/(2n\sigma_1^{4})]$ leads to the priors
$$
\begin{array}{c}
\mu|\sigma^2,H_0 \sim N(0,g_0\sigma^2), \qquad
\sigma^2|H_0 \sim LN(0,2g_0), 
\\ [1ex]
\mu_0|\sigma_1^2,H_1 \sim N(0,g_1(n/n_0)\sigma_0^2), \qquad 
\mu_1|\sigma_1^2,H_1 \sim N(0,g_1(n/n_1)\sigma_1^2), 
\\ [1ex]
\sigma_0^2|H_1 \sim LN(0,2(n/n_0)g_1), \quad \mbox{and} \quad
\sigma_1^2|H_1 \sim LN(0,2(n/n_1)g_1).
\end{array}
$$

\noindent Setting $g_j = h^{1/d_j}$, $j=0,1$ results in
$$
\begin{array}{l}
\ds p(\vx|H_0) \\
\begin{array}{cl}

\ = & \ds \int_0^\infty
\frac{1}{\sqrt{2\pi}^n}\exp\left[
- \tfrac{n+2}{2}\ln(\sigma^2)
- \tfrac{n\widehat{\sigma}_{h^{1/2}}^2}{2\sigma^2}
- \tfrac{\ln(4\pi h)}{2}
- \tfrac{(\ln \sigma^2)^2}{4h}
- \tfrac{1}{2}\ln\left(n + \tfrac{1}{\sqrt{h}}\right)
\right] d\sigma^2
\\ [2ex]
\  \stackrel{h\Rightarrow\infty}{=} 
& \ds \int_0^\infty
\frac{1}{\sqrt{2\pi}^n}\exp\left[
- \tfrac{n+2}{2}\ln(\sigma^2)
- \tfrac{n\widehat{\sigma}^2}{2\sigma^2}
- \tfrac{\ln(4\pi)}{2}
- \tfrac{1}{2}\ln(n)
\right] d\sigma^2
\\ [2ex]


\  = & \ds
\exp\left[ \ln p(\vx|\widehat{\vtheta}_0)
+ \xi\left( \tfrac{n}{2} \right)
- \tfrac{1}{2}\ln(2)
- \tfrac{1}{2}\ln(n)
\right],
\end{array}
\end{array}
$$

\noindent where
$\xi(x) = \ln\Gamma(x) + x - x\ln(x) - (1/2)\ln(2\pi)$
and $\ln p(\vx|\widehat{\vtheta}_0) = - (n/2)\ln\left( 2\pi\widehat{\sigma}^2 \right)
- n/2$ is the log-likelihood of the null model
evaluated at its MLE $\widehat{\vtheta}_0$.
Similarly, 
$p(\vx_0|H_1,h)$ is 
$$
\begin{array}{rcl}
\ds p(\vx_0|H_1)
%
%
%

& = & \ds \int_0^\infty  \tfrac{1}{\sqrt{2\pi}^{n_0}}
\exp\Big[  
- \left(\tfrac{n_0}{2} + 1\right)\ln( \sigma_0^2)
- \tfrac{1}{2\sigma_0^2}
\left\{ \| \vx_0\|^2 
- \tfrac{(n_0\overline{x}_0)^2}{n_0 + (n_0/n)h^{-1/4}} 
\right\} 
\\
&  & \ds \qquad 
- \tfrac{1}{2}\ln(n + h^{-1/4})
- \tfrac{1}{2}\ln(2\pi(2n/n_0)h^{1/2})
- \tfrac{(\ln \sigma_0^2)^2}{4(n/n_0) h^{1/4}}
\Big]   d\sigma_0^2.

\end{array}
$$

\noindent 
By construction the $\ln(h^{1/2})$ term 
cancels in the numerator and denominator of the Bayes factor.
The integrand is a monotonic function of $h$ 
(apart from the $\ln(h^{1/2})$ term which cancels) and
has a well defined limit as $h\to\infty$. Taking
$h\to\infty$ the above expression for  $p(\vx_0|H_1,h)$  
simplifies to 
$p(\vx_0|H_1)
%
= \exp\left[
\ell(\widehat{\mu}_0,\widehat{\sigma}_0^2)
+ \xi\left( \tfrac{n_0}{2} \right)
- \tfrac{1}{2}\ln(2n^2/n_0)
\right],
$
%
%
where $\ell(\widehat{\mu}_0,\widehat{\sigma}_0^2)
= - (n_0/2)\ln( 2\pi\widehat{\sigma}_{0}^2) - n_0/2$. Combining with a 
similarly obtained expression for $p(\vx_1|H_1)$ we obtain
$$
\begin{array}{rl}
\lambda_{\mbox{\scriptsize Bayes}} 
& \ds = \lambda_{\mbox{\scriptsize LRT}} 
- 3\ln(n)
- 2\xi(n/2)
+ 2\xi(n_0/2)
+ 2\xi(n_1/2)
+ \ln(n_0 n_1/2),
\end{array} 
$$
\noindent 
where
$\lambda_{\mbox{\scriptsize LRT}}
= n\ln(\widehat{\sigma}^2)
- n_0\ln(\widehat{\sigma}_{0}^2)
- n_1\ln(\widehat{\sigma}_{1}^2)$, the estimators $\widehat{\sigma}^2 = n^{-1}\|\vx - \overline{x}\vone\|^2$,
$\widehat{\sigma}_{0}^2 = n_0^{-1}\|\vx_0 - \overline{x}_0\vone\|^2$,
and
$\widehat{\sigma}_{1}^2 = n_1^{-1}\|\vx_1 - \overline{x}_1\vone\|^2$
are the MLEs for the variance parameters.
Stirling's asymptotic expansion
of $\ln\Gamma(z)$ for large
$z$ is
$\ln \Gamma (z) = z\ln(z)-z -(1/2)\ln(z) + (1/2)\ln(2\pi) +O(z^{-1})$ \citep[see for example Equation 6.1.37 of][]{Abramowitz1972}. Hence,
$\xi(x) = - \tfrac{1}{2}\ln(x) + O(x^{-1})$.
Using this   $\lambda_{\mbox{\scriptsize Bayes}}$ simplifies to
$$
\lambda_{\mbox{\scriptsize Bayes}}  = \lambda_{\mbox{\scriptsize LRT}} - 2\ln(n) + O(n_0^{-1} + n_1^{-1}).
$$

\noindent
Note that the coefficient of $\ln(n)$ is $d_1 - d_0 = 2$
which is the corresponding degrees of freedom of the corresponding LRT.

The parameter posteriors are given by:
$$
\begin{array}{c}
\mu|\vx,H_0 \sim t_n(\overline{x},n^{-1}\widehat{\sigma}^2),
\quad
\sigma_0^2|\vx,H_0 \sim \mbox{IG}\left(\tfrac{n}{2},\tfrac{n}{2}\widehat{\sigma}^2\right),
\quad 
\mu_0|\vx,H_1 \sim t_{n_0}(\overline{x}_0,n_0^{-1}\widehat{\sigma}_{0}^2)
\\ [1ex]
\sigma_0^2|\vx,H_1 \sim \mbox{IG}\left(\tfrac{n_0}{2},\tfrac{n_0}{2}\widehat{\sigma}_{0}^2\right), \quad 
\mu_1|\vx,H_1 \sim t_{n_1}(\overline{x}_1,n_1^{-1}\widehat{\sigma}_{1}^2), \quad \mbox{and} \quad
\sigma_1^2|\vx,H_1 \sim \mbox{IG}\left(\tfrac{n_1}{2},\tfrac{n_1}{2}\widehat{\sigma}_{1}^2\right).
\end{array}
$$


It is important to note that all of the
constant terms have cancelled from the asymptotic
approximation for $\lambda_{\mbox{\scriptsize Bayes}}$.
This has been achieved by incorporating the
$(n/n_0)$ and $(n/n_1)$ factor in the priors for $\mu_0$ and $\mu_1$,
and the $(2n/n_0)$ and $(2n/n_1)$ factors in the priors
for $\sigma_0^2$ and $\sigma_1^2$. Without these factors,
cancellation of $O(1)$ and larger terms in the expression
$\lambda_{\mbox{\scriptsize Bayes}}$ would not occur.

We conducted a short numerical experiment to compare our Bayesian test with the LRT. We simulated
$10^6$ datasets where $n_0 = n_1 = 50$ with the true parameter values
\begin{enumerate}
\item[(i)] $\mu_0 = 0$, $\sigma_0 = \sigma_1 = 1$ and $\mu_1 \in \{0, 0.25, 0.5, 1\}$; or

\item[(ii)] $\mu_0 = 0$, $\mu_1 =0$,  $\sigma_0 = 1$ and $\sigma_1 \in \{0, 1.25, 1.5, 2.5 \}$.
\end{enumerate}

\begin{figure}[h!]
	\centering
	\includegraphics[width=0.88\textwidth]{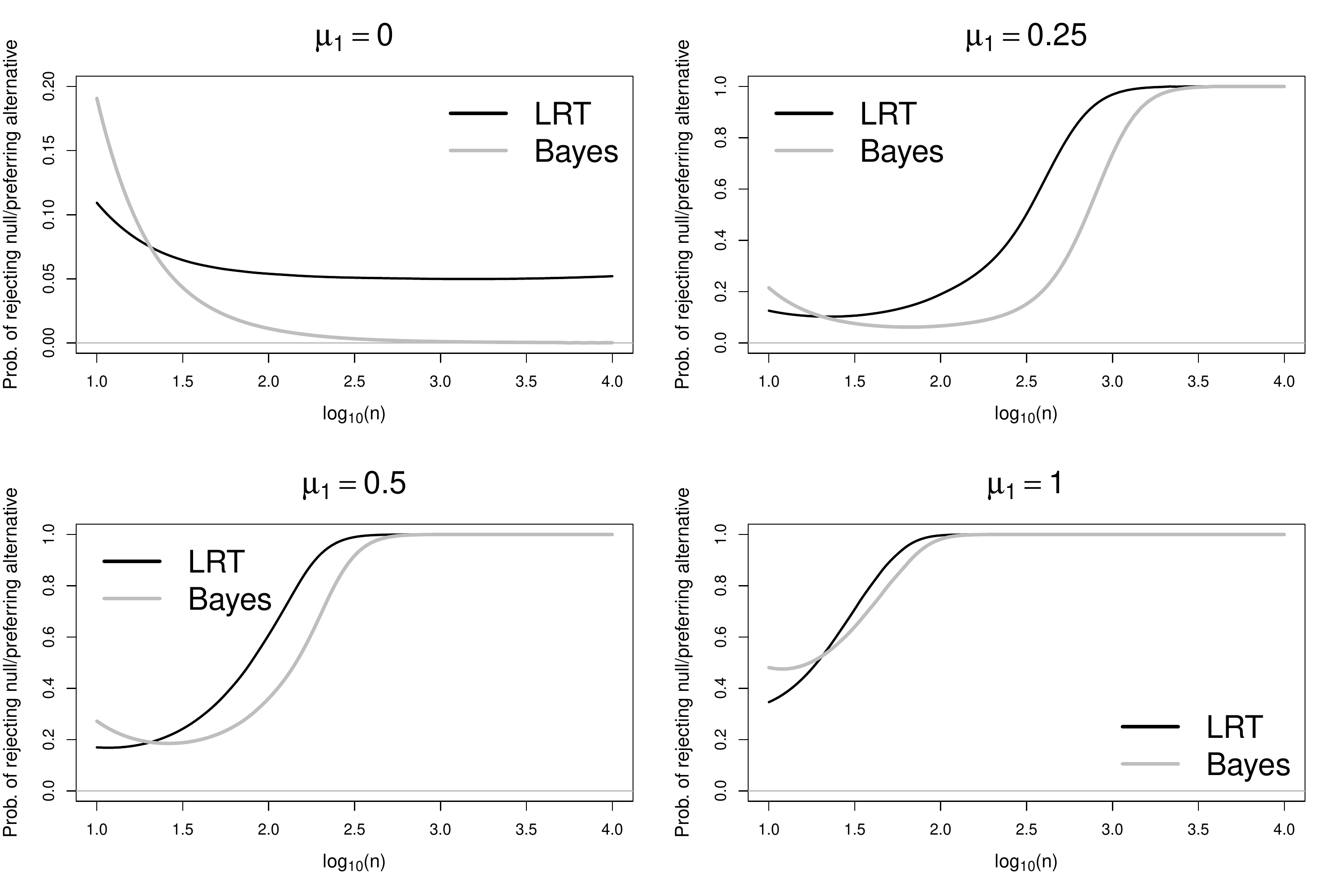}
	
	\includegraphics[width=0.88\textwidth]{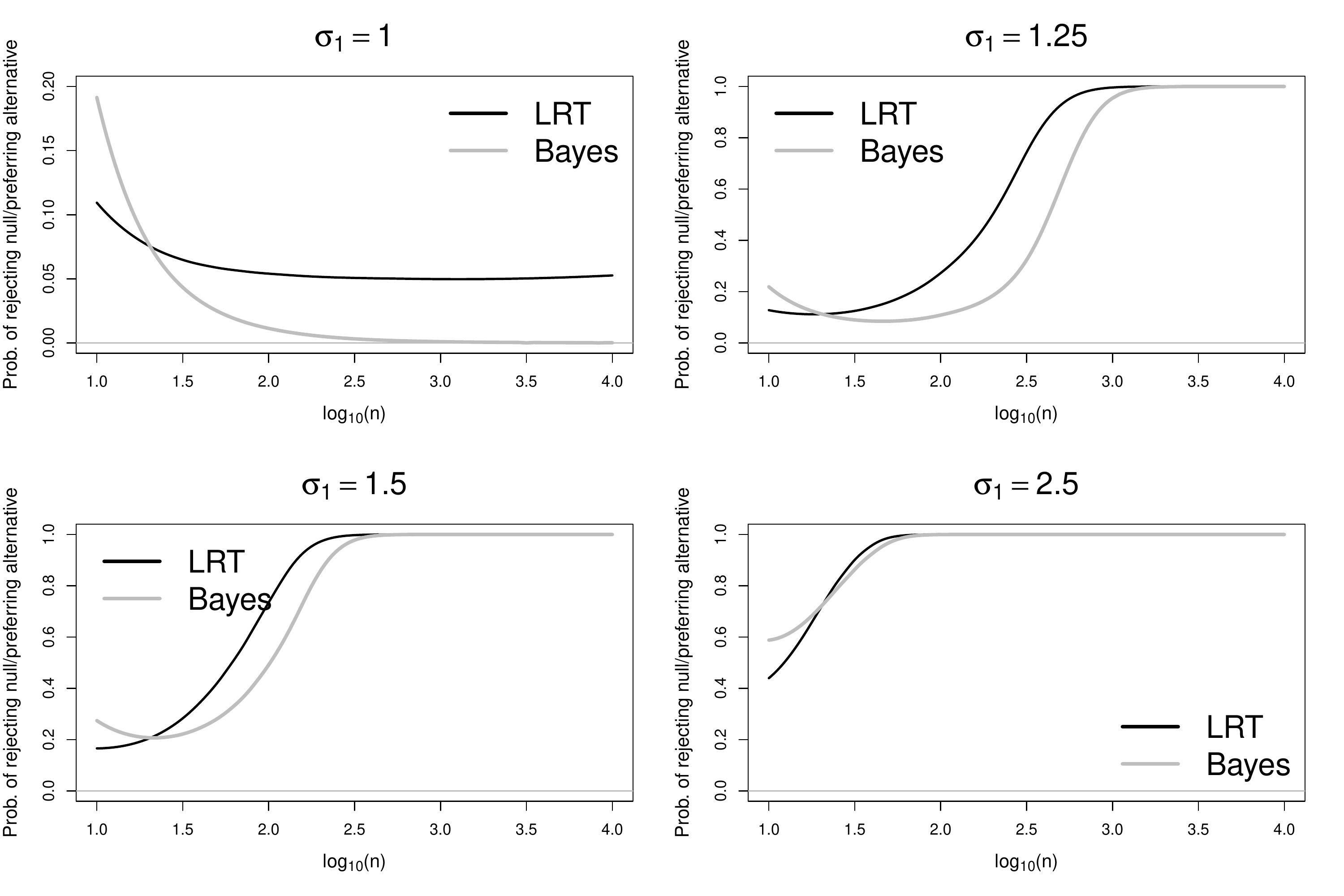}
	\caption{The empirical probabilities of rejecting the null
		(in the LRT case) or preferring the alternative (in the Bayesian
		test) when simulating two normal populations with $n_0 = n_1 = 50$, and (i)
		$\mu_0 = 0$, $\sigma_0 = \sigma_1 = 1$ and $\mu_1 \in \{0, 0.25, 0.5, 1\}$ (left four panels);
		or (ii) $\mu_0 = 0$, $\mu_1 =0$,  $\sigma_0 = 1$ and $\sigma_1 \in \{0, 1.25, 1.5, 2.5 \}$
		(right four panels).}
	\label{fig:03}
\end{figure}

\noindent 
The empirical probabilities of rejecting the null
(in the LRT case) or preferring the alternative (in the Bayesian
test) are illustrated in Figure \ref{fig:03}. Note that
under $H_0$ the type I error approaches zero as $n\to \infty$
for the Bayesian test,
and under $H_1$ the type II error approaches zero as $n\to \infty$
for both the Bayesian and LRT tests. When $H_0$ is true the
LRT has a fixed 5\% type I error.

\subsection{Linear models}
\label{sec:Example2}

We will now consider hypothesis testing for linear models.
Consider the base model
$$
\ds \vy|\alpha,\vbeta,\sigma^2
\sim N(
\alpha\vone + \mX\vbeta,\sigma^2\mI),
$$

\noindent where
$\vy$ is a response vector of length $n$, $\vbeta$ is a coefficient vector of length $p$, $\sigma^2$ is a positive scalar,
$\mX$ is a full-rank $n$ by $p$ matrix of covariates, and $\mI$ is the identity matrix of appropriate dimension. In order to simplify some calculations we will transform $\vy$ and $\mX$
so that $\vy$ and the columns of $\mX$ are standardized, i.e., $\overline{y} = 0$, $\|\vy\|^2 = \vy^T\vy = n$, $\mX_j^T\vone = 0$, and $\|\mX_j\|^2 = n$ where $\mX_j$ is the $j$th column of $\mX$. Let $\vgamma$ be a binary vector of length $p$, and
let $\mX_\vgamma$ be the submatrix $\mX$
comprised from the columns of $\mX$ whose
corresponding elements of $\vgamma$ are non-zero.
Consider the hypothesis  test
\begin{equation}\label{eq:hypothesesLM}
H_0 \colon \vgamma = \vgamma_0
\qquad \mbox{versus} \qquad
H_1 \colon \vgamma = \vgamma_1,
\end{equation}

\noindent where $\vgamma_0$ and $\vgamma_1$ denote the models under the null and alternative hypotheses respectively with $0\le|\vgamma_0|\le|\vgamma_1|$.

To simplify exposition for this example we will only use cake
priors for $\alpha$ and $\vbeta_\vgamma$. Since
$\sigma^2$ is a common parameter across all parameters
we can use the typical improper (Jeffreys) priors for $\sigma^2$ given by
$p(\sigma^2) \propto (\sigma^2)^{-1} \bI(\sigma^2>0)$.
This choice has been formally justified in \cite{Berger1998}.
Cake priors can be used for all parameters for this example,
but the working out is lengthy and unnecessarily obfuscates the exposition.
Using 
$\mP(\alpha,\vbeta_\vgamma) = \widetilde{\sI}(\alpha,\vbeta_\vgamma) = \mbox{diag}(\sigma^{-2},\sigma^{-2}\mX_\vgamma^T\mX_\vgamma/n)$
for a particular model $\vgamma$ leads to
\begin{equation}\label{eq:proirs2}
\ds \alpha|\sigma^2,g \sim N(0,g\sigma^2), \quad \mbox{and} \quad
\ds \vbeta_\vgamma|\sigma^2,g \sim N\left(\vzero,g\sigma^2\left(\tfrac{1}{n}\mX_\vgamma^T\mX_\vgamma\right)^{-1} \right).
\end{equation}

\noindent Further, we use
$p(\vbeta_{-\vgamma}) = \prod_{j\colon\gamma_j=0} \delta(\beta_j;0)$ 
where $\delta(x;a)$ is the Dirac delta function with location $a$.
The prior on $\vbeta_\vgamma$ is simply the Zellner $g$-prior \citep{Zellner1986}
 where the prior covariance is scaled by a factor of $n$.
The prior on $\vbeta_\vgamma$ combined with the prior on $\vbeta_{-\vgamma}$ is a spike and slab prior
for $\vbeta$.
%

Marginalizing over $\alpha$, $\vbeta$ and $\sigma^2$
for a particular model $\vgamma$ we obtain after simplification
$$
p(\vy|\vgamma,g)
= \frac{\Gamma(n/2)}{(n\pi)^{n/2}}\exp\Big[
- \tfrac{1+|\vgamma|}{2}\ln(g)
- \tfrac{1+|\vgamma|}{2}\ln(n + g^{-1})
- \tfrac{n}{2}\ln\left( 1 - \tfrac{g}{1+g} R_{\vgamma}^2 \right)
\Big],
$$

\noindent where 
$R_\vgamma^2$ 
 is the usual R-squared statistic for model $\vgamma$.
This is equivalent to the $g=h^{1/d_h}$ of Section \ref{sec:Cake}  
to using the hyperpriors
$p(g|\vgamma_j) = \delta(g; h^{1/(1 + |\vgamma_j|)})$, $j=0,1$.
After marginalizing over $g$
the Bayes factor
as a function of $h$ simplifies to
$$
\begin{array}{rl}
\ds \mbox{BF}_{01}(h)
& \ds =
\exp\Big[
- \tfrac{n}{2}\ln\left( 1 - \tfrac{h^{1/(1+|\vgamma_0|)}}{1 + h^{1/(1+|\vgamma_0|)}} R_{\vgamma_0}^2 \right)
+ \tfrac{n}{2}\ln\left( 1 - \tfrac{h^{1/(1+|\vgamma_1|)}}{1 + h^{1/(1+|\vgamma_1|)}} R_{\vgamma_1}^2 \right) \\ [2ex]
& \ds \qquad \qquad \qquad
- \tfrac{1+|\vgamma_0|}{2}\ln\left(n + h^{-1/(1+|\vgamma_0|)} \right)
+ \tfrac{1+|\vgamma_1|}{2}\ln\left(n + h^{-1/(1+|\vgamma_1|)} \right)
\Big].
\end{array}
$$

\noindent
Taking $h\to \infty$
we use the fact that 
$1 - R_\vgamma^2 = \widehat{\sigma}_{\vgamma}^2$
(where $\widehat{\sigma}_{\vgamma}^2$ is the MLE for $\sigma^2$ under the
model $\vgamma$)
to obtain
$$
\begin{array}{rl} 
\ds \lambda_{\mbox{\scriptsize Bayes}}
& \ds = 
\left[ - n\ln( \widehat{\sigma}_{\vgamma_1}^2) - |\vgamma_1| \ln(n) \right]
- \left[  - n\ln( \widehat{\sigma}_{\vgamma_0}^2) - |\vgamma_0| \ln(n) \right]
\\ [1ex]
& \ds =  \mbox{BIC}_{\vgamma_0} - \mbox{BIC}_{\vgamma_1} 
= \lambda_{\mbox{\scriptsize LRT}} - \nu\ln(n),
\end{array} 
$$

\noindent where
$\ds \mbox{BIC}_{\vgamma}
=
n\ln(2\pi\widehat{\sigma}_{\vgamma}^2)
- n + |\vgamma|\ln(n)
= -2 \ln p(\vy|
\widehat{\alpha}_{\vgamma},
\widehat{\vbeta}_{\vgamma}, \widehat{\sigma}_{\vgamma}^2) + |\vgamma|\ln(n),
$
and, $\widehat{\alpha}_{\vgamma}$ and $\widehat{\vbeta}_{\vgamma}$ are the MLEs
for $\alpha$ and $\vbeta$ under model $\vgamma$, $\lambda_{\mbox{\scriptsize LRT}} = n\ln( \widehat{\sigma}_{\vgamma_0}^2 )
- n\ln( \widehat{\sigma}_{\vgamma_1}^2 )$
is the LRT statistic corresponding to the
hypotheses (\ref{eq:hypothesesLM}) and
$\nu = |\vgamma_1| - |\vgamma_0|$.
Hence, for these models and prior structures the Bayesian
test statistic is simply the difference between two BIC values.

%
%

Note that as $h\to\infty$ the parameter posteriors become
$$\alpha|\vy,\vgamma \sim t_n(0,\widehat{\sigma}_{\vgamma}^2/n), \quad
\vbeta_{\vgamma}|\vy,\vgamma \sim t_n( \widehat{\vbeta}_{\vgamma}, \widehat{\sigma}_{\vgamma}^2 \left(\mX_\vgamma^T\mX_\vgamma  \right)^{-1} ),
\quad \mbox{and} \quad  
\sigma^2|\vy,\vgamma \sim \mbox{IG}\left( \tfrac{n}{2}, \tfrac{n}{2}\widehat{\sigma}_{\vgamma}^2 \right),
$$

\noindent 
where $\widehat{\vbeta}_{\vgamma}$
and $\widehat{\sigma}_{\vgamma}^2$ are the
MLEs corresponding to model $\vgamma$.

We will not provide any numerical examples due to the close relationship
between our Bayes factors and the BIC, 
and the fact that almost every paper ever written on model selection for linear models uses the BIC
in its comparisons. 
We direct the interested reader to any
of the papers in the discussion below all of which make comparisons with the BIC as a model selection criteria.

There are four main differences between the priors used here and the
priors that have been used in the literature for linear models.
The first such difference is the choice of prior on $\alpha$  
which the typical prior is to use the Jeffreys prior $p(\alpha) \propto 1$ which was
advocated in \cite{Berger1998}.
If we were to use this prior and were
only to use cake priors for $\vbeta_\vgamma$ then
$p(g_0|\vgamma_0;h) = \delta(g_0;h^{1/|\vgamma_0|})$ instead of
$p(g_0|\vgamma_0;h) = \delta(g_0;h^{1/(1+|\vgamma_0|)})$. The consequence of
this would be that the null model (where $\vgamma = \vzero$) would become problematic to calculate.

The second difference is in the choice of prior for $\vbeta$.
Most Bayesian approaches to model selection for linear models use the Zellner $g$-prior
where
\begin{equation}\label{eq:alternativeBeta}
\vbeta_\vgamma|\sigma^2,g \sim N(\vzero,g\sigma^2(\mX_\vgamma^T\mX_\vgamma)^{-1}),
\end{equation}

\noindent instead of the prior for $\vbeta_\vgamma$ in (\ref{eq:proirs2}).
This difference is subtle.

\cite{Bayarri2012} advocate the priors for $\vbeta$ 
should remain proper and not degenerate to a point mass.
If we were to treat $\mX_\vgamma$ as random then under mild conditions
$\mX_\vgamma^T\mX_\vgamma/n \to \bE(\mX_\vgamma^T\mX_\vgamma)$
almost surely suggesting that our prior for $\vbeta_\vgamma$
does not degenerate to a point mass. Lastly
(\ref{eq:proirs2}) and (\ref{eq:alternativeBeta})
simplify marginal likelihoods
since terms involving determinants cancel, and have the added advantage that
they do
not depend on the unit of measurements  of the covariates.

The third difference is the choice of prior on $g$.
As stated in the introduction, \cite{Liang2008} argues for a hyperprior to be assigned to $g$.
\cite{Liang2008} considers the hyper $g$-prior; 
the hyper $g/n$-prior;
and the Zellner-Siow prior (equivalent to a particular inverse-gamma on $g$)  
\citep{Zellner1980b}. \cite{Maruyama2011} use a different prior to
(\ref{eq:proirs2}) or (\ref{eq:alternativeBeta})
and a beta-prime prior with specially chosen prior hyperparameter values.
All of these choices, apart from \cite{Maruyama2011},  
either no closed form expression for the marginal likelihood 
exists, or such an expression is in terms of a Gauss hypergeometric function
which is numerically difficult to evaluate \citep{Pearson2017} so that approximation is required. 

Model selection consistency is another desirable criteria of
\cite{Bayarri2012}.
The authors corresponding $g/n$, Zellner-Siow, beta-prior and robust
priors are model selection consistent for all possible models.
Our prior specification results in a null based Bayes factor
is a simple function of the BIC and so
achieves model selection consistency for iid data \citep[under
some
additional mild assumptions,][]{Yang2005}.
See also Section \ref{sec:comparingTwoModels}.

\subsection{Handling zero parameters in the null model}
\label{sec:zeroCase}

\noindent
Let us now return to Lindley's example posed by \cite{Lindley1957}
described in Section \ref{sec:LindlyBardletParadox}.
In order to apply the methodology of Section \ref{sec:Cake}
the null model needs to have a non-zero number of parameters.
We provide the following novel artificial construct to handle this case in order
to augment the problems so that both hypotheses have a non-zero
number of parameters.
\begin{enumerate}
	\item Introduce a second sample of hypothetical data, say $\vz$.
	
	\item Modify the null and alternative hypotheses by adding a clause
	that the hypothetical data has the same distribution under the null and
	alternative hypotheses.
	
	\item Apply the methodology of Section \ref{sec:Cake} to
	the augmented problem.
\end{enumerate}

\noindent

\noindent In order to illustrate this approach suppose
we have a second sample of hypothetical data
$\vz = (z_1,\ldots,z_n)^T$
and consider the augmented  hypotheses
$$
\begin{array}{l}
H_0\colon x_1,\ldots,x_n \sim N(\mu_0,\sigma^2) \qquad \mbox{and} \qquad
z_1,\ldots,z_n|\widetilde{\mu} \sim N(\widetilde{\mu},\sigma^2) \qquad \mbox{versus} 
\\ [1ex]
H_1\colon x_1,\ldots,x_n|\mu \sim N(\mu,\sigma^2) \qquad \mbox{and} \qquad
z_1,\ldots,z_n|\widetilde{\mu} \sim N(\widetilde{\mu},\sigma^2),
\end{array}
$$

\noindent where $\mu_0$ and $\sigma^2$ have  known fixed values,
and $\widetilde{\mu}$ is an artificial mean parameter corresponding
to the sample $\vz$.
This is a modification of the original hypotheses (\ref{eq:hypZtest})
has the same logical implication as
the hypotheses (\ref{eq:hypZtest})
for the observed sample $\vx$ since the hypothetical data has the same
hypothetical models under the null and alternative hypotheses.
For the augmented problem we have
$\vtheta_0 = \widetilde{\mu}$
with $d_0 = 1$, and $\vtheta_1 = (\mu,\widetilde{\mu})^T$
and $d_1 = 2$ so that we have avoided the problem of dividing by
zero.
The cake priors become
$\widetilde{\mu}|H_0 \sim N(0,g_0\sigma^2)$, 
 $\mu|H_1 \sim N(0,g_1\sigma^2)
$ and $
\widetilde{\mu}|H_1 \sim N(0,g_1\sigma^2).
$
For $g_0$ and $g_1$ we use
$g_0 =h$ and
$g_1 = h^{1/2}$.
Then
$$
\begin{array}{rl}
\ln p(\vx|H_0)

& \ds =
- \tfrac{n}{2}\ln(2\pi\sigma^2)
- \tfrac{\| \vx - \mu_0 \vone \|^2}{2\sigma^2},

\\ [2ex]

\ln p(\vz|H_0)
& \ds =
- \tfrac{n}{2}\ln(2\pi\sigma^2)
- \tfrac{1}{2\sigma^2} \left[ \|\vz\|^2 - \tfrac{(n\overline{z})^2}{n + h^{-1}} \right]
- \tfrac{1}{2}\ln(h)
- \tfrac{1}{2}\ln(n+ h^{-1}),

\\ [2ex]

\ln p(\vx|H_1)
& \ds =
- \tfrac{n}{2}\ln(2\pi\sigma^2)
- \tfrac{1}{2\sigma^2} \left[ \|\vx\|^2 - \tfrac{(n\overline{x})^2}{n + h^{-1/2}} \right]
- \tfrac{1}{2}\ln(h^{1/2})
- \tfrac{1}{2}\ln\left(n+ \tfrac{1}{\sqrt{h}}\right), \ \mbox{and}

\\ [2ex]

\ln p(\vz|H_1)
& \ds =
- \tfrac{n}{2}\ln(2\pi\sigma^2)
- \tfrac{1}{2\sigma^2} \left[ \|\vz\|^2 - \tfrac{(n\overline{z})^2}{n + h^{-1/2}} \right]
- \tfrac{1}{2}\ln(h^{1/2})
- \tfrac{1}{2}\ln\left(n+ \tfrac{1}{\sqrt{h}}\right).
\end{array}
$$

\noindent The Bayes factor in the limit as $h\to\infty$ is
$$
\lambda_{\mbox{\scriptsize Bayes}}
= \lim_{h\to\infty}-2\ln\left[  \frac{p(\vx|H_0)p(\vz|H_0) }{p(\vx|H_1)p(\vz|H_1)  } \right] = \lambda_{\mbox{\scriptsize LRT}} - \ln(n),
$$

\noindent
where
$\lambda_{\mbox{\scriptsize LRT}} = \sigma^{-2} (\| \vx - \mu_0 \vone \|^2
- \| \vx - \widehat{\mu} \vone \|^2)$ is the likelihood ratio
test statistic corresponding to the hypothesis (\ref{eq:hypZtest}).
We conducted a small simulation study identical to the
simulation study in Section \ref{sec:Example1} with the exception
that $\sigma^2$ was treated as known. The resulting figure and
interpretation was nearly identical to that
in Section \ref{sec:Example1} (not shown).


\section{Theory}
\label{sec:theory}

In all of the examples in Section \ref{sec:Cake} the
quantity $\lambda_{\mbox{\scriptsize Bayes}}$ can be placed into the form
(\ref{eq:BayesLRT}).
We will now consider the asymptotic properties of hypothesis tests
based on this form.
\cite{Shao2003}
developed theory regarding the asymptotic properties of hypothesis tests. We will adopt his
notation and definitions here.
Let $\mX = (X_1, ...,X_n)^T$ be a random sample from $\sP = \{ \, p_i( \, \cdot \,) \colon i=1,\ldots,n \, \}$.
%
The type I and type II errors are defined by
$\alpha_T(\sP) = \bP(\, T(\mX) = 1)$ 
when $\sP\in\sP_0$   
and
$1 - \alpha_T(\sP) = \bP(\,T(\mX) = 0)$ 
 when $\sP\in\sP_1$   
%
\noindent respectively.
Fix the level of significance $\alpha$ such that $\sup_{\sP\in\sP_0} \{ \alpha_T({\sP}) \} \leq \alpha$.
We will now suppose that $T_n(\mX) \equiv T(\mX)$
and consider scenarios where $n$ diverges.
In our ensuing discussion we use the following definitions.

\medskip
\noindent
{\bf Definitions from 2.13 of \cite{Shao2003}:}
\begin{enumerate}
\item[(i)] If $\ds \lim_{n\to\infty} \sup_{\sP\in\sP_0} \{ \alpha_{T_n}(\sP) \} \le \alpha$  then $\alpha$ is an asymptotic significance level of $T_n$.

\item[(ii)] If
$\ds \lim_{n\to\infty} \sup_{\sP\in\sP_0} \{ \alpha_{T_n}(\sP) \}$ exists,
then it is called the limiting size of $T_n$.

\item [(iii)] The sequence of tests $T_n$ is called consistent if and only if the type II error probability converges to $0$, i.e.,
$\lim_{n\to\infty}[1 - \alpha_{T_n}(\sP)] = 0$, for any
$\sP\in\sP_1$.

\item [(iv)] The sequence of tests $T_n$ is called Chernoff-consistent if and only if
$T_n$ is consistent and the type I error probability converges to 0, i.e., $\lim_{n\to\infty} \{ \alpha_{T_n}(\sP)\} = 0$, for any
$\sP\in \sP_0$. Furthermore, $T_n$ is called strongly Chernoff-consistent if and only if $T_n$ is consistent
and the limiting size of $T_n$ is $0$.
\end{enumerate}


\noindent 
We note that any reasonable test which is consistent where the level $\alpha$ is
controllable,  can be made Chernoff-consistent by letting $\alpha_n\equiv\alpha \to 0$
as $n\to\infty$.

Wilks Theorem \citep{Wilks1938} tells
us that, assuming the data was generated
under the null distribution, under appropriate
regularity conditions (including that the hypotheses
are nested) that
$\lambda_{\mbox{\scriptsize LRT}}$ converges  to $\chi_{\nu}^2$ in distribution so that $\lambda_{\mbox{\scriptsize LRT}} = O_p(1)$. 
A detailed exposition on the characterization the asymptotic distribution of the LRT statistic under quite general conditions, including when $H_0$ and/or $H_1$ is misspecified, and whether the hypotheses are nested
or non-nested, can be found in \cite{vuong89}. Below we summarize the most pertinent
results.

\subsection{Asymptotic properties of the likelihood ratio test statistic} 

Suppose \(X_1,X_2,\ldots,X_n\) are independent random variables
from  $\{ \, p_{0i}(\,\cdot\,) \colon i=1,\ldots,n \, \}$ and that we have a parametric model
\(\mathcal{P}=\{\, p_i(\,\cdot\, | \vtheta) \colon i=1,\ldots,n,\
\vtheta\in\Theta \, \}\) (which may or may not include the true distribution(s) \(\{\,p_{0i}(\,\cdot\,)\, \}\)).
Define the log-likelihood as
$\ell(\vtheta) =\sum_{i=1}^n \ln p_i(X_i|\vtheta)$,
the MLE and ``pseudo-true'' value of $\vtheta$  as
$$
\ds \widehat{\vtheta} = \arg\max_{\theta} \, \{ \, \ell(\vtheta) \, \}
\qquad \mbox{and} \qquad 
\vtheta^* = \arg\max_{\vtheta} \ \left[ \bE\left\{ n^{-1}\ell(\vtheta) \right\} \right],
$$

\noindent respectively (assuming both are well-defined). 
Assume also that
$\bE\left[ n^{-1}\ell(\vtheta^*) \right]\to C,$ 
for some finite \(C>0\). 
Under the ``pseudo-true'' value $\vtheta^*$ 
	the resultant distribution
is $\{ \, p_i(\,\cdot\, | \vtheta^*) \, \}$ which is the ``best'' distribution in the sense that it results in the smallest Kullback-Leibler (KL) divergence
with respect to the true distribution over all distributions for the model.
We assume that these are unique.

If \(\mathcal{P}\) is suitably regular, 
with \(\Theta\) a nice subset of
\(d\)-dimensional Euclidean space,
then certain derivatives exist and various statements
can be made: the Euclidean norm of
\(\widehat{\vtheta}-\vtheta^*\) is \(O_p(n^{-1/2})\), 
and writing \( \nabla \ell(\vtheta)\) and \( \nabla^2 \ell(\vtheta)\) for the first and
second order partial derivatives we may expand \(\ell(\vtheta^*)\)
about \(\widehat{\vtheta}\) to get
$$
\begin{array}{rl}
\ell(\vtheta^*)
& = 
\ell(\widehat{\vtheta})
+ (\vtheta^* - \widehat{\vtheta})^T \nabla \ell(\widehat{\vtheta}) 
+ \tfrac{1}{2}(\vtheta^*-\widehat{\vtheta})^T 
\left[ \nabla^2 \ell(\widetilde{\vtheta}) \right] 
(\vtheta^*-\widehat{\vtheta}) \\
& = \ell(\widehat{\vtheta})
- \tfrac{1}{2}\Big[ n^{1/2}(\widehat{\vtheta}-\vtheta^*) \Big]^T  
\bE\Big[ - \tfrac{1}{n} \nabla^2 \ell(\vtheta^*)\Big] \Big[ n^{1/2}(\widehat{\vtheta}-\vtheta^*) \Big] + o_p(1),
\end{array} 
$$

\noindent since \( \nabla \ell(\widehat{\vtheta})\equiv0\), 
for some \(\widetilde{\vtheta}\) between \(\vtheta^*\) and \(\widehat{\vtheta}\); also
the quadratic form is \(O_p(1)\). Thus, we may decompose the maximised
log-likelihood into the following terms:
$$
\ds \ell(\widehat{\vtheta}) =  
\bE \left[ \ell(\vtheta^*) \right]
 +  \left[ \ell(\vtheta^*) - \bE \left\{ \ell(\vtheta^*) \right\} \right]  
 + O_p(1).
$$

\noindent 
The first term on the right hand side is asymptotic to \(nC\); the second term is a random sum
of \(n\) terms
with expectation zero, and so under further regularity conditions is
\(O_p(n^{1/2})\). We refer to these three terms as the ``\(O_p(n)\)'',
``\(O_p(n^{1/2})\)'' and ``\(O_p(1)\)'' terms respectively (from left to right).
Precise regularity conditions guaranteeing all of this nice behaviour
may be found in 
\cite{vuong89}; see also Chapter 5 of \cite{vaart98}. They are all
easily satisfied in all of our examples.

Finally, we note that the first term can be written as:
$$
\begin{array}{rl} 
\ds \bE \left[ \ell(\vtheta) \right] 
& \ds = \bE \left[ \sum_{i=1}^n \ln\left\{
\frac{p_i(X_i|\vtheta)}{p_{0i}(X_i)} \right\}   \right] + \bE\left[ \sum_{i=1}^n \ln p_{0i}(X_i) \right]
\\
& \ds \widehat{=} - \mbox{KL}( \, p_{0i} \, || \, p_i  \, ) 
+ \bE\left[ \sum_{i=1}^n \ln p_{0i}(X_i) \right]
\end{array} 
$$

\noindent where the first term corresponds to the negative 
KL-divergence between $p_{0i}$ and $p_i$,
and the second term is related to the entropy of $p_{0i}$. We see from the above
equation that maximizing $\bE \left[ \ell(\vtheta) \right]$  with respect to $\vtheta$ is equivalent to minimizing 
the KL-divergence between $\prod_{i=1}^n p_{0i}(\,\cdot\,)$ and 
$\prod_{i=1}^n p_i(\,\cdot\,|\vtheta)$.

\subsection{Comparing models}
\label{sec:comparingTwoModels}

Suppose now that we have two competing models \(\mathcal{P}_j = \{ \, p_{ij}(\,\cdot\, | \vtheta_j,H_j) \colon i=1,\ldots,n,\
\vtheta_j\in\Theta \, \}, j = 0, 1\),  with corresponding log-likelihoods \(\ell_0(\,\cdot\,)\) and
\(\ell_1(\,\cdot\,)\), MLEs \(\widehat{\vtheta}_0 \)
and \(\widehat{\vtheta}_1\), and pseudo-true values 
\(\vtheta_0^*\) and 
\(\vtheta_1^*\). For convenience we will write \(\widehat{\ell}_j \equiv 
\ell_j(\widehat{\vtheta}_j) \). 
Suppose also that for \(j=0,1\), we have 
$\bE \left[ n^{-1}\ell_j(\vtheta_j^*) \right] \to \ell_j^*$.

There are two main cases:
\begin{enumerate}
	\item If $\ell_1^* > \ell_0^*$ then the model under
$H_1$ has a smaller KL-divergence from the
true distribution than the model under $H_0$, and
	we have immediately
	\begin{align*}
	n^{-1}\lambda_{\mbox{\scriptsize LRT}}(\mX) =  n^{-1} (
	\widehat \ell_1 -
	\widehat \ell_0
	)\stackrel P\to \ell_1^* - \ell_0^* >0\,, 
	\end{align*}
	that is the LRT statistic
	$\lambda_{\mbox{\scriptsize LRT}}(\mX) = 2(\widehat \ell_1 - \widehat \ell_0)$ 
	is of order \(n\) in probability.
	
	\item If \(\ell_1^* = \ell_0^*\) then immediately we
	see that the two ``\(O_p(n)\)'' terms in the difference between the
	log-likelihoods would, at least asymptotically, cancel out and that
	the ``\(O_p(n^{1/2})\)'' terms would ``dominate''. However, in many
	practical examples \emph{the ``\(O_p(n^{1/2})\)'' terms also cancel out}, in
	which case
	\begin{align*}
	\lambda_{\mbox{\scriptsize LRT}}(\mX) = 2(\widehat \ell_1 - \widehat \ell_0)
	=O_p(1)\,.
	\end{align*}
	
	\citep[see in particular Theorem 3.3 of][]{vuong89}. This occurs when the ``best'' member of each model yields the \emph{same
		distribution}, that is when
	$\{\,p_{i0}(\,\cdot\,|\vtheta_0^*,H_0) \, \}=
	\{\, p_{i1}(\,\cdot\,|\vtheta_1^*,H_1) \, \}$.
	The parametrisations may be completely different, but nonetheless the
	distributions
	corresponding to the pseudo-true parameter values are identical.
	This occurs when
	the two models have some overlap, i.e., \(\mathcal{P}_0 \cap
	\mathcal{P}_1 \) is non-empty \emph{as a subset of all possible sets of
		distributions} \(  \{ \, p_i(\,\cdot \, ) \, \}\). 
	In such a case, the models may or may not be nested, and may or may
	not be correctly specified; however the ``best'' distribution in both is
	the same (and is part of the overlap).
\end{enumerate}

\noindent 
We also have the following consequences.
\begin{itemize}
	
	\item {\bf Lindley's paradox:}   	If $H_1$ is correct (and $H_0$ is not) then the above theory implies $\lambda_{\mbox{\scriptsize LRT}}(\mX)$ is $O_p(n)$ and
	a test of the form (\ref{eq:BayesLRT}) is consistent. Since Lindley's paradox occur with asymptotic probability $p(\chi_{\nu,\alpha}^2 < \lambda_{\mbox{\scriptsize LRT}} < \nu \ln n)$,
		it occurs with vanishingly small probability
		as $n\to \infty$.
	
	\item {\bf Chernoff-consistency:} 
	If $H_0$ is true then $\lambda_{\mbox{\scriptsize LRT}}(\mX)$ is $O_p(1)$
	and $\bP(T(\mX) = 1) \to 0$ as $n\to\infty$. If $H_1$ is true then 
	$\lambda_{\mbox{\scriptsize LRT}}(\mX)$ is $O_p(n)$
	and $\bP(T(\mX) = 0) \to 0$ as $n\to\infty$. Hence, a Bayesian test of the form
	(\ref{eq:BayesLRT}) is Chernoff consistent.

	\item {\bf Model selection consistency:} 
	Now suppose that we have $M$ competing hypotheses $H_j$, $j=1,\ldots,M$. For each model we have
	$\bE \left[ n^{-1}\ell_j(\vtheta_j^*) \right] \to \ell_j^*$ for some constants $\ell_j^*$.
	We will call a hypothesis $H_j$ correct if $j\in\sC$ where
	$$
	\ds \sC = \Big\{ \, j \, \colon  \, \ell_j^* = \max_{k=1,\ldots,M}  \ell_k^* \, \Big\}.
	$$
	
	\noindent The hypotheses $H_j$ such that $j\in\sC$ correspond to correct 
	models in the sense that all such models are closest in terms of their
	KL-divergence to the data generating distribution. 
	Using (\ref{eq:BayesLRT}) to compare models not in $\sC$ with models in $\sC$ leads to $\lambda_{\mbox{\scriptsize LRT}}(\mX)$ is $O_p(n)$ and the test preferring the model in $\sC$. Comparing any two models in
	$\sC$ will asymptotically prefer the model with the smallest size.
\end{itemize}	

\noindent 
When cake priors become arbitrarily diffuse to the point of becoming improper and a further problem occurs.
We now discuss such problems.
 
\section{Arbitrary constants}
\label{sec:ArbitraryConstants}

We now return to the issue of  improper priors in the context of Bayes factors discussed in the introduction.
When using improper priors, consider the
conditional  density $p(\vx|\vtheta_i,H_j)$   where
$p(\vtheta|H_j) \propto f_i(\vtheta_j)   =D_j f_j(\vtheta_j)
$
for some $D_j>0$
such that $\int f_j(\vtheta_j) d\vtheta_j$ does not exist for $j=0,1$. Then
\begin{equation}\label{eq:BayesFactorWithImproperPriors}
\ds B_{01} = \frac{D_0}{D_1} \times \frac{\int p(\vx|\vtheta,H_0) f_0(\vtheta_0)d\vtheta_0}{\int p(\vx|\vtheta,H_1) f_1(\vtheta_1)d\vtheta_1}.
\end{equation}

\noindent This Bayes factor is problematic since it depends on two arbitrary constants $D_0$
and $D_1$.
\begin{itemize}
	\item {\bf Problem III}: When using improper
	priors either the null or the alternative model can be made
	to be preferred by artificially changing $D_0$ or $D_1$
	to suit the a priori preferred conclusion.
\end{itemize}

\noindent In the limit as as $h\to\infty$ cake priors become improper.
Suppose that instead of using $g_j = h_j^{1/d_j}$ we use
$g_j= (D_jh)^{1/d_j}$, where $c_j>0$ are arbitrary constants. This implies
$$
\lambda_{\mbox{\scriptsize Bayes}} = \lambda_{\mbox{\scriptsize LRT}} - \nu\ln(n) + \Delta + o_p(1)
$$

\noindent where $\Delta = \ln(D_1/D_0)$
is a controllable constant determined by how
the $g_j$ parameters diverge. Thus, we have not avoided the problem of arbitrary 
constants in our test.

Based on the theory developed in Section \ref{sec:theory}, if $\Delta = O(1)$
then the corresponding test will have all of the properties discussed in
Section (\ref{sec:comparingTwoModels}) where the level of the test requires
adjustment.
For small sample sizes the value of $\Delta$ trades-off the probabilities of type I
and type II errors.
We do not believe that the model nor the data can
make the choice of $\Delta$ value on behalf of the analyst. Furthermore, any criteria used to
select $\Delta$ is making the choice in trade-off of relative probabilities of type I
and type II errors whether implicitly or explicitly.
Implicitly cake priors are choosing $\Delta = 0$. 

This choice might be preferred for the following reasons:
\begin{enumerate}
	\item 
This choice leads to a parsimonious expression for $\lambda_{\mbox{\scriptsize Bayes}}$ (all $O(1)$ terms cancel); 

\item When the number of parameters in the null and alternative hypotheses are
		equal ($\nu=0$) using $\Delta=0$ leads
		to preferring the model with the larger 
		likelihood. If $\Delta\ne 0$ this is always the case. 
		If the analyst desired to explicitly favour 
		ether model when $\nu=0$ then 
		the prior odds should be altered to achieve this; and 
		
\item 
Choosing different values of $D_0$ and $D_1$ would be inconsistent with
		similar choices made in the literature. For example,
		in the linear models example one might also use $p(\alpha^2|H_j) = D_j$ or $p(\sigma^2|H_j) = D_j (\sigma^2)^{-1} I(\sigma^2 > 0)$. No papers in the model selection
		literature, to our knowledge, chose different constants $D_j$
		for each model under consideration. This
		suggests choosing $D_0 = D_1$ is reasonable.
\end{enumerate}		

Lastly, we note that if we select $\Delta = \nu\ln(n) - \chi_{\nu,\alpha}^2$, where 
$\chi_{\nu,\alpha}^2$ denotes the upper quantile function of the chi-squared distribution
with degrees of freedom parameter $\nu$, that Bayesian tests can be made to mimic
the frequentists LRT when the type I error is controlled to have level $\alpha$.

\section{Interpretation}
\label{sec:interpretation}

For Bayesian test using cake priors and the LRT test functions are (approximately)
\begin{equation}\label{eq:both} 
\ds T_{\mbox{\scriptsize Bayes}}(\vx) = I(\lambda_{\mbox{\scriptsize LRT}}(\vx) > \nu\ln(n))
\qquad \mbox{and} \qquad 
T_{\mbox{\scriptsize LRT}}(\vx) = I(\lambda_{\mbox{\scriptsize LRT}}(\vx) > \chi_{\nu,\alpha}^2)
\end{equation}

\noindent respectively. This also allows us to obtain a one-to-one mapping between
$p$-values and values of $\lambda_{\mbox{\scriptsize Bayes}}(\vx)$.
It can be shown
\begin{equation}\label{eq:chisquardBounds} 
\ds \nu + 2\ln(1/\alpha) - 5/2 \le \chi_{\nu,\alpha}^2 \le \nu + 2\ln(1/\alpha) + 2\sqrt{\nu\ln(1/\alpha)},
\end{equation}

\noindent where the lower bound only holds for $\alpha \le 0.17$ and $\nu\ge 2$
\citep{Laurent2000,Inglot2010}. 
%
%
%
For our Bayesian tests the cut-off value, $\nu\ln(n)$ grows with $n$ implying
that the level of the test decays with $n$.
A consequence is that our Bayesian tests offer some protection against
a sequential analysis (where samples are collected and statistical significance
checked sequentially).
Suppose that $\alpha = (n/e)^{\nu/2}$ then (\ref{eq:chisquardBounds}) becomes
$\ds \nu\ln(n) - 5/2 \le \chi_{\nu,\alpha}^2 \le \nu\ln(n) + \sqrt{2\nu^2\ln(n)}$
so that $\chi_{\nu,\alpha}^2 \asymp \nu\ln(n)$,
and the frequentist and Bayesian testing procedures are roughly equivalent.

Given that (\ref{eq:both}) we can directly compare 
that the Bayes factors have taken the interpretation offered by \cite{Kass1995}
in Table \ref{tab:bayesFactorInterpretation}
appears to have the short-coming of not taking
into account the degrees of freedom $\nu$ of the test, nor the sample
size $n$. Consider Table \ref{tab:BaysFactors} which directly compares
$p$-values with their corresponding $\lambda_{\mbox{\scriptsize Bayes}}$
for given
$n$ and $\nu$. If one were to use a LRT instead where
$\lambda_{\mbox{\scriptsize Bayes}}(\vx)$ takes the threshold values in Table
\ref{tab:bayesFactorInterpretation}, i.e., $\lambda_{\mbox{\scriptsize Bayes}}(\vx)\in \{ 0, 2,6,10\}$,
$\nu$ ranges from 1 to 5, and $n\in \{ 50, 10^2, 10^3 \}$.
Note that every $p$-value is
smaller than 5\% suggesting that Bayesian tests are typically more conservative at
preferring the alternative than classical tests reject the null at the usual 5\% cut-off.
Further, going from $\lambda_{\mbox{\scriptsize Bayes}}(\vx) = 2$ to
$\lambda_{\mbox{\scriptsize Bayes}}(\vx) = 6$ and from $\lambda_{\mbox{\scriptsize Bayes}} = 6$ to
$\lambda_{\mbox{\scriptsize Bayes}}(\vx) = 10$ roughly translates
to a 5 to 10 fold reduction in the corresponding $p$-value.
Thus, the anticipated potential short-coming of Table \ref{tab:bayesFactorInterpretation}
not depending on $\nu$ or $n$ does not pan out, at least for the values of $n$ and
$\nu$ considered in Table \ref{tab:BaysFactors}. We believe 
Table \ref{tab:BaysFactors} is a reasonable interpretation of strength of evidence
for Bayesian tests.  


\begin{table}[h]
	\centering
	\small{ 	
		\begin{tabular}{|c|r|rrr|}	
			\hline
			& & \multicolumn{3}{c|}{$p$-values}
			\\
			\hline
			$\nu$ &	$\lambda_{\mbox{\tiny Bayes}}$ &			    $n=50$ & $n=10^2$ & $n=10^3$  \\
			\hline
			1&	0 &            4.79E-02 & 3.18E-02 & 8.58E-04 \\
			1&	2 &			   1.50E-02 & 1.02E-02 & 2.84E-03  \\
			1&	6 &			   1.64E-03 & 1.13E-03 & 3.27E-04 \\
			1&	10 &	       1.92E-04 & 1.33E-04 & 3.92E-05  \\
			\hline
			2&	0&                 2.00E-02 & 1.00E-02 & 1.00E-03 \\
			2&	2	&			   7.36E-03 & 3.68E-03 & 3.68E-04  \\
			2&	6		&			   9.96E-04 & 4.98E-04 & 4.98E-05  \\
			2&	10		&			   1.35E-04 & 6.74E-05 & 6.74E-06 \\
			\hline	
			3&	0       &              8.34E-02 & 3.16E-03 & 1.20E-04 \\
			3&	2		&		   3.29E-03 & 1.24E-03 & 4.61E-05 \\
			3&	6		&			   4.99E-04 & 1.85E-04 & 6.73E-06  \\
			3&	10		&			   7.40E-05 & 2.73E-05 & 9.72E-07  \\
			\hline		
			4&	0       &             3.53E-03 & 1.02E-03 & 1.48E-05 \\
			4&	2		&			  1.45E-03 & 4.12E-04 & 5.82E-06  \\
			4&	6		&			  2.35E-04 & 6.58E-05 & 8.87E-07  \\
			4&	10		&			   3.73E-05 & 1.02E-05 & 1.34E-07  \\
			\hline
		\end{tabular}\begin{tabular}{l|rrr|}
			\hline		
			& \multicolumn{3}{c|}{$\lambda_{\mbox{\scriptsize Bayes}}$}
			\\
			\hline
			$p$-value &			 $n=50$ & $n=10^2$ & $n=10^3$  \\
			\hline
			0.05        &			   {$-$0.1} & {$-$0.8} & {$-$3.1}  \\
			0.01       &			    2.7 &  2.0 & {$-$0.3} \\
			0.001      &			    {6.9} &  {6.2} &  3.9  \\
			0.0001     &              {11.2} & {10.5} & {8.2} \\
			\hline
			0.05	    &			   {$-$1.8} & {$-$3.2} & {$-$7.8}  \\
			0.01		&			   1.4 &  0.0 & {$-$4.6}  \\
			0.001		&			   {6.0} &  4.6 & 0.0 \\
			0.0001      &              {10.6}  & {9.2} &  4.6 \\
			\hline	
			0.05		&			   {$-$3.9} & {$-$6.0} & {$-$12.9} \\
			0.01		&			   {$-$0.4} & {$-$2.5} & {$-$9.4}  \\
			0.001		&			   4.5  & 2.5 & {$-$4.5}  \\
			0.0001      &              {9.4} & {7.3} & 0.4 \\
			\hline		
			0.05		&			    {$-$6.2} & {$-$8.9} & {$-$18.1}  \\
			0.01		&			    {$-$2.4} & {$-$5.1} & {$-$14.4}  \\
			0.001		&			     2.8 &  0.0 &  {$-$9.2}  \\
			0.0001      &              {7.9}  & 5.1 & {$-$4.1} \\ 
			\hline
		\end{tabular}
	} 
	\caption{The third to sixth columns are LRT based $p$-value for 
		the $\lambda_{\mbox{\scriptsize Bayes}} \in \{ 0, 2, 6, 10\}$ specified in
		the second column for different values of $n$ and $\nu$.
		The eighth to 11th colums are values of $\lambda_{\mbox{\scriptsize Bayes}}$
		for the $p$-value specified in the seventh column for different values of $n$ and $\nu$. }
	\label{tab:BaysFactors}
\end{table}

\section{Conclusion}
\label{sec:discussion}

We have introduced a new class of priors we call cake priors having a number of desirable properties.
Cake priors   can be made arbitrarily diffuse without
the Bayes factor favouring the null or alternative
hypotheses as the prior becomes increasingly diffuse.
In the limit, at least for the examples we consider here, Bayes factors take the form of penalized likelihood ratio statistics, one of the most thoroughly
understood quantities in Statistics. Due to their close
link with Jeffreys priors, cake priors are parametrization invariant. The resulting Bayesian test avoids the need to specify a $p$-value cut-off and
are asymptotically Chernoff-consistent. With a slight
modification, an arbitrary but controllable constant
  can be used to bridge the gap between Bayesian
tests and likelihood ratio tests. 
Unlike approaches that split the
dataset up into parts or use imaginary data,
cake priors are transparent, uncomplicated, and easily implementable.
Finally, Bayesian tests using cake priors providing some protection against sequential testing being more conservative as the sample size grows.
We believe all of the above properties should make cake priors the
default choice when performing Bayesian hypothesis tests
for hypothesis consisting of
a simple point null against a composite alternative for parametric models with
iid data.

\bibliographystyle{elsarticle-harv}
\bibliography{CakeBib}

\end{document}